\documentclass[12pt,reqno,final]{amsart}

\usepackage{amsmath}
\usepackage{graphicx}
\usepackage{algorithmic}
\usepackage{algorithm}
\usepackage{amssymb}
\usepackage[tight,footnotesize]{subfigure}
\usepackage[scaled]{helvet} 
\usepackage{fourier}
\usepackage{fullpage}

\DeclareGraphicsExtensions{.eps}
\usepackage[usenames]{color}
\def\x{{\mathbf x}}

\def\P {{\mathbf{P}}}         
\def\Q {{\mathbf{Q}}}         
\def\CF{{\mathcal{S}}}
\def\CFI{{\mathcal{S}}}
\def\CFII{{\tilde{\mathcal{S}}}}
\def\psc{{\zeta}}
\def\pscI{{\zeta}}
\def\pscII{{\tilde{\zeta}}}
\def\CI{{C}}
\def\CII{{\tilde{C}}}

\newtheorem{thm}{Theorem}[section]
\newtheorem{cor}[thm]{Corollary}
\newtheorem{lem}[thm]{Lemma}
\newtheorem{prop}[thm]{Proposition}
\newtheorem{defn}[thm]{Definition}

\newcommand{\ceil}[1]{\left\lceil#1\right\rceil}
\newcommand{\norm}[1]{\left\Vert#1\right\Vert}
\newcommand{\abs}[1]{\left\vert#1\right\vert}
\newcommand\vect[1]{{\bf#1}}
\newcommand\matr[1]{{\bf#1}}

\newcommand\alphabf{{\boldsymbol{\alpha}}}
\newcommand{\argmin}{\operatornamewithlimits{argmin}}
\newcommand{\argmax}{\operatornamewithlimits{argmax}}

\newcommand{\Real}{\mathbb R}
\newcommand\RR[1]{\mathbb{R}^{#1}}

\DeclareMathOperator{\supp}{supp}
\DeclareMathOperator{\ext}{ext}

\DeclareMathOperator{\range}{range}

\begin{document}

\title{Greedy Signal Space Methods for incoherence and beyond} 

\author{Raja Giryes and Deanna Needell}


\begin{abstract}
Compressive sampling (CoSa) has provided many methods for signal recovery of signals compressible with respect to an orthonormal basis.  However, modern applications have sparked the emergence of approaches for signals not sparse in an orthonormal basis but in some arbitrary, perhaps highly overcomplete, dictionary.  Recently, several ``signal-space'' greedy methods have been proposed to address signal recovery in this setting.  However, such methods inherently rely on the existence of fast and accurate projections which allow one to identify the most relevant atoms in a dictionary for any given signal, up to a very strict accuracy.  When the dictionary is highly overcomplete, no such projections are currently known; the requirements on such projections do not even hold for incoherent or well-behaved dictionaries.   In this work, we provide an alternate analysis for signal space greedy methods which enforce assumptions on these projections which hold in several settings including those when the dictionary is incoherent or structurally coherent.  These results align more closely with traditional results in the standard CoSa literature and improve upon previous work in the signal space setting.
\end{abstract}

\maketitle


\section{Introduction}
\label{sec:intro}

In many signal and image processing applications we encounter the following problem: recovering an original signal $\vect{x}\in \RR{d}$ from a set of noisy measurements
\begin{eqnarray}
\label{eq:measurement}
\vect{y} = \matr{M}\vect{x}+ \vect{e},
\end{eqnarray}
where $\vect{M}\in \RR{m\times d}$ is a known linear operator and $\vect{e} \in \RR{d}$ is additive bounded noise, i.e. $\norm{\vect{e}}_2^2 \le \varepsilon^2$.  In many cases such as those in Compressive Sampling (CoSa)~\cite{DSPweb}, we have $m\ll d$ and thus \eqref{eq:measurement} has infinitely many solutions.  
To make the problem well-posed we rely on additional priors for the signal $\vect{x}$, such as \textit{sparsity}.

The sparsity assumption provides two main models, termed the
\textit{synthesis} and \textit{analysis} models \cite{elad07Analysis}.
The synthesis model, which has received great attention in the past decade,
assumes that $\vect{x}$ has a $k$-sparse representation $\alphabf$ under a given dictionary $\matr{D} \in \RR{d \times n}$ \cite{Bruckstein09From}.
In other words, there exists a vector $\alphabf \in \RR{n}$ such that $\vect{x} = \matr{D}\alphabf$ and $\norm{\alphabf}_0 \le k$,
where $\norm{\alphabf}_0 = |\supp(\alphabf)|$ denotes the $\ell_0$ pseudo-norm.
Under the synthesis model assumption we can recover $\vect{x} = \matr{D}\alphabf$ by solving
\begin{eqnarray}
\label{eq:l0_synthesis}
\argmin_\alphabf \norm{\alphabf}_0 & s.t. & \norm{\vect{y} - \matr{M}\matr{D}\alphabf}_2 \le \varepsilon.
\end{eqnarray}

Since solving \eqref{eq:l0_synthesis} is an NP-complete problem in general~\cite{Davis97Adaptive}, approximation techniques are required for recovering $\vect{x}$. One strategy uses relaxation, replacing the $\ell_0$ with the $\ell_1$ norm, resulting in the $\ell_1$-synthesis problem
\begin{eqnarray}
\label{eq:l1_synthesis}
\hat{\alphabf}_{\ell_1} = \argmin_{\alphabf} \norm{\alphabf}_1 & s.t. & \norm{\vect{y} - \matr{M}\matr{D}\alphabf}_2 \le \varepsilon.
\end{eqnarray}
The study of these types of synthesis programs has largely relied on properties like the \textit{Restricted Isometry Property} (RIP)~\cite{Candes05Decoding}, which states that
$$
(1-\delta_{k})\norm{\vect{x}}^2 \leq \norm{\matr{M}\vect{x}}^2 \leq (1+\delta_k)\norm{\vect{x}}^2 \quad\text{for all $k$-sparse $\vect{x}$},
$$
for some small enough constant $\delta_k < 1$.

If the matrix $\matr{D}$ is unitary and the vector $\vect{x}$ has a $k$-sparse representation $\alphabf$, then when $\matr{M}$ satisfies the RIP with $\delta_{2k} < \delta_{\ell_1}$, the program~\eqref{eq:l1_synthesis} accurately recovers the signal,
\begin{eqnarray}
\label{eq:l1_rec_error}
\norm{\hat{\vect{x}}_{\ell_1} - \vect{x}}_2 \le C_{\ell_1}\varepsilon,
\end{eqnarray}
where $\hat{\vect{x}}_{\ell_1} = \matr{D}\hat{\alphabf}_{\ell_1}$, 
$C_{\ell_1}$ is a constant greater than $\sqrt{2}$ and $\delta_{\ell_1}$ ($\simeq 0.4652$) is a constant \cite{Candes06Stable,Candes06Near, foucart10Sparse}.  This result also implies perfect recovery in the absence of noise.
It was extended also for incoherent redundant dictionaries \cite{Rauhut08Compressed}.

An alternative aproach to approximating \eqref{eq:l0_synthesis} is to use a greedy strategy. 
Recently introduced methods that use this strategy are the CoSaMP \cite{Needell09CoSaMP},
  IHT \cite{Blumensath09Iterative}, and HTP \cite{Foucart11Hard} methods.  Greedy methods iteratively identify elements of the support
 of the signal, and once identified, use a simple least-squares to recover the signal.
These methods were shown to have guarantees in the form of
\eqref{eq:l1_rec_error} under the assumption of the RIP. 
However, such results hold only when $\matr{D}$ is orthonormal, and do not hold for general dictionaries $\matr{D}$.
Recently, the greedy approaches have been adapted to this setting.  For example, the \textit{Signal Space CoSaMP} method \cite{Davenport13Signal} adapts CoSaMP to the setting of arbitrary dictionaries.  A slight modification\footnote{Here we use two separate support selection schemes, whereas the original Signal Space CoSaMP method uses one.} of this method is shown in Algorithm~\ref{alg:SSCoSaMP}.  In the algorithm, the subscript $T$ denotes the restriction to elements (columns) indexed in $T$.
The function $\CF_k(\vect{y})$ returns the support of the best $k$-sparse representation of $\vect{y}$ in the dictionary $\matr{D}$, and $\matr{P}_T$ denotes the projection onto that support.

\begin{algorithm}[t]
\caption{Signal Space CoSaMP (SSCoSaMP)} \label{alg:SSCoSaMP}
\begin{algorithmic}[l]

\REQUIRE $k, \matr{M}, \matr{D}, \vect{y}, a$ where $\vect{y} = \matr{M}\vect{x}
+ \vect{e}$, $k$ is the sparsity of $\vect{x}$ under $\matr{D}$ and $\vect{e}$ is
the additive noise.  $\CFI_{k}$ and $\CFII_{ak}$ is a pair of near optimal projection schemes.

\ENSURE $\hat{\vect{x}}$: $k$-sparse approximation of
$\vect{x}$.

\STATE Initialize the support $T^0 =\emptyset$, the residual $\vect{y}_r^0 = \vect{y}$ and set $t = 0$.

\WHILE{halting criterion is not satisfied}

\STATE $t = t + 1$.

\STATE Find new support elements: $T_\Delta=\CFII_{ak}(\matr{M}^*\vect{y}^{t - 1}_r)$.

\STATE Update the support: $\tilde{T}^t = T^{t -1} \cup
T_\Delta$.

\STATE Compute the representation: $\vect{x}_p = \matr{D}(\matr{MD}_{\tilde{T}^{t}})^{\dag}\vect{y} = \matr{D}\left(\argmin_{\tilde\alphabf}\norm{\vect{y} -\matr{MD}\tilde\alphabf}_2^2 \text{ s.t. } \tilde\alphabf_{(\tilde{T}^t)^C}=0\right)$.

\STATE Shrink support: $T^t = \CFI_{k}(\vect{x}_p)$.

\STATE Calculate new representation: ${\vect{x}}^t = \P_{T^t}\vect{x}_p$.

\STATE Update the residual:
$\vect{y}_r^t = \vect{y} - \matr{M}{\vect{x}}^t$.

\ENDWHILE
\STATE Form final solution $\hat{\vect{x}} = {\vect{x}}^t$.
\end{algorithmic}
\end{algorithm}

In~\cite{Davenport13Signal}, the authors analyze this CoSaMP variant under the assumption of the \textit{$\matr{D}$-RIP} \cite{Candes11Compressed}, which states\footnote{By abuse of notation we denote both the RIP and the $\matr{D}$-RIP constants by $\delta_k$. It will be clear from the context to which one we refer at each point in the article.}
\begin{equation}\label{def:D_RIP}
(1-\delta_{k})\norm{\matr{D}\alphabf}^2 \leq \norm{\matr{M}\matr{D}\alphabf}^2 \leq (1+\delta_k)\norm{\matr{D}\alphabf}^2 \quad\text{for all $k$-sparse $\alphabf$}.
\end{equation}
They prove that under this assumption, if one has access to projections $\CF_k$ which satisfy
\begin{equation}\label{dnw_reqs}
\norm{\CF_k(\vect{z}) - \CF_k^{opt}(\vect{z})}_2 \leq \min\left(c_1\norm{\CF_k^{opt}(\vect{z})}_2, c_2\norm{\vect{z} - \CF_k^{opt}(\vect{z})}_2 \right),
\end{equation}
then the method accurately recovers the $k$-sparse signal, as in~\eqref{eq:l1_rec_error}.  Here, we write $\CF_k^{opt}$ to denote the optimal projection.  Although the other results for greedy methods in this setting also rely on similar assumptions~\cite{Blumensath11Sampling,Giryes13Greedy}, it remains an open problem whether such projections can be obtained.  In this paper, we address this issue by analyzing the two projections in the method separately, and using an alternative theoretical analysis.  This analysis allows us to weaken the requirement on the projections.  This new requirement also shows that when $\matr{D}$ is incoherent, traditional compressed sensing algorithms can be used for these projections.  Of course, the interesting case is when the dictionary $\matr{D}$ is not incoherent at all.

\subsection{Our contribution} In this paper we present a variant of SSCoSaMP and develop theoretical guarantees for it.  We provide similar guarantees to CoSaMP for incoherent dictionaries and show how these are extended for coherent ones. 

As is evident by Algorithm~\ref{alg:SSCoSaMP}, as in the case of other greedy methods, we need access to a projection which, given a general vector, finds the closest (in the $\ell_2$ sense) $k$-sparse vector.

In the representation case (when $\matr{D} = \matr{I}$), simple hard thresholding gives the desired result.
However, in the signal space we need to solve
\begin{eqnarray}
\label{eq:optimal_sparse_projection}
\CF^*_{k}(\vect{z}) = \argmin_{\abs{T}\le k} \norm{\vect{z} - \P_T\vect{z}}_2^2.
\end{eqnarray}
This problem seems to be NP-hard in general, as is the case in the analysis framework \cite{Gribonval13Projection}, so an approximation is needed.
For this we introduce the notion of a pair of\textit{ near-optimal projection}, which extends the definition in \cite{Giryes13Greedy} and is similar to the one in \cite{Davenport13Signal}.
\begin{defn}
\label{def:C_optimal_proj}
A pair of procedures ${\CFI}_{\pscI k}$ and ${\CFII}_{{\pscII} k}$ implies a pair of near-optimal projections $\P_{{\CFI}_{\pscI k}(\cdot)}$ and $\P_{{\CFII}_{{\pscII} k}(\cdot)}$ with constants $\CI_k$ and ${\CII}_k$ if for any $\vect{z} \in \Real^d$, $\abs{\CFI_{\pscI k}(\vect{z})} \le \pscI k$, with $\pscI \ge 1$, $\abs{{\CFII}_{{\pscII} k}(\vect{z})} \le {\pscII} k$, with ${\pscII} \ge 1$, and
\begin{eqnarray}
\label{eq:C_optimal_proj}
 \norm{\vect{z}-\P_{\CFI_{\pscI k}(\vect{z})}\vect{z}}_2^2 \le \CI_k\norm{\vect{z} - \P_{\CF^*_k(\vect{z})}\vect{z}}_2^2 \quad\text{as well as}\quad
 \norm{\P_{{\CFII}_{{\pscII} k}(\vect{z})}\vect{z}}_2^2 \ge {\CII}_k\norm{\P_{\CF^*_k(\vect{z})}\vect{z}}_2^2,
\end{eqnarray}
where $\P_{\CF^*_k}$ denotes the optimal projection as in~\eqref{eq:optimal_sparse_projection}.
\end{defn}

Our main result can now be summarized as follows.

\begin{thm}
\label{thm:general_bound}
Let $\vect{y} = \matr{M}\x +\vect{e}$, where $\matr{M}$ satisfies the $\matr{D}$-RIP~\eqref{def:D_RIP} with a constant $\delta_{(3\pscI+1)k}$ ($\pscI\ge 1$), $\x$ be a vector with a $k$-sparse representation under $\matr{D}$ and $\vect{e}$ is a vector of additive noise.  Suppose that $\CFI_{\pscI k}$ and $\CFII_{2\pscI k}$ are a pair of near optimal projections (as in Definition~\ref{def:C_optimal_proj}) with constants $\CI_k$ and ${\CII}_{2k}$.
Apply SSCoSaMP (with $a = 2$) and let ${\vect{x}}^t$ denote the approximation after $t$ iterations.  If
$\delta_{(3\pscI+1) k}< \epsilon^2_{\CI_k,\CII_{2k},\gamma}$ and 
 \begin{eqnarray}
\label{eq:C_k_tilda_C_2k_cond}
\left(1+\sqrt{\CI_k} \right)^2\left( 1-\frac{{\CII}_{2k}}{(1+\gamma)^2} \right)<1,
\end{eqnarray}
then after a constant number of iterations $t^*$ it holds that
\begin{eqnarray}
\label{eq:general_bound}
&& \hspace{-0.5in} \norm{{\vect{x}}^{t^*} -\vect{x}}_2 \le   \eta_0\norm{\vect{e}}_2,
\end{eqnarray}
where $\gamma$ is an arbitrary constant, and
$\eta_0$ is a constant depending on $\delta_{(3\pscI+1)k}$, $\CI_k$, ${\CII}_{2k}$ and $\gamma$.  
The constant $\epsilon_{\CI_k,\CII_{2k},\gamma}$ is dependent on the values $\CI_k$, $\CII_{2k}$, 
and $\gamma$; it is the solution to a quadratic equation involving these parameters which is 
 greater than zero if and 
only if \eqref{eq:C_k_tilda_C_2k_cond} holds.
\end{thm}

{\bfseries Remark.}  Note that we use $a=2$ as in the traditional CoSaMP method for our analysis, but like in the traditional method, Algorithm~\ref{alg:SSCoSaMP} provides a template and other choices of $a$ ($a\geq 1)$ can certainly be used.  Similarly, a large value of $\pscI$ allows~\eqref{eq:C_k_tilda_C_2k_cond} to be satisfied with smaller constants, but makes the RIP restriction $\delta_{(3\pscI+1) k}< \epsilon^2_{\CI_k,{\CII}_{2k},\gamma}$ much harder to be satisfied.  This is an important tradeoff, and different projections will optimize this tradeoff in different ways.  To allow for a wide range of projections to be used, we leave $\pscI$ as a free parameter.

Unlike previous results in the signal space setting, the requirement~\eqref{eq:C_k_tilda_C_2k_cond} on the near-optimal projections holds in many common compressed sensing settings such as those when the dictionary $\matr{D}$ is incoherent or satisfies the RIP.  In those settings, classical recovery methods may be utilized for the projections.  We thus offer an improvement over existing signal space analyses which enforce requirements on the projections which do not even hold when the dictionary is highly incoherent.

\subsection{Organization}
In Section~\ref{sec:notation} we present the notation we use in the work
and some preliminaries.
We present the proof of our main result, Theorem~\ref{thm:general_bound}, in Section~\ref{sec:guarantees}.  In Section~\ref{sec:near_opt_proj} we examine some important settings to which our results apply.
Section~\ref{sec:conc} discusses related works and concludes the paper.


\section{Notations and Preliminaries}
\label{sec:notation}

We use the following notation in our work.  We write $\norm{\cdot}_2$ for the Euclidean ($\ell_2$) norm of a vector, and $\norm{\cdot}$ for the spectral ($\ell_2 \rightarrow \ell_2$) norm of a matrix.  We denote the identity matrix by $\matr{I}_d = \matr{I} \in \RR{d \times d}$. Given support set $T$, $\matr{D}_{T}$ is the sub-matrix of $\matr{D}$ whose columns are indexed by $T$.   We write $\P_T = \matr{D}_T\matr{D}_T^\dag$ as the orthogonal projection onto $\range(\matr{D}_T)$ and
$\Q_T = \matr{I}_d - \P_T$ as the orthogonal projection onto the orthogonal
complement of $\range(\matr{D}_T)$.

Recall the $\matr{D}$-RIP defined in~\eqref{def:D_RIP}, which enforces that the measurement matrix $\matr{M}$ preserves the geometry of signals sparse with respect to $\matr{D}$.  The $\matr{D}$-RIP, like the standard RIP, inherits the following useful properties.
The first follows immediately from the definition and thus appears without a proof.

\begin{cor}
\label{cor:MP_RIP_norm}
If $\;\matr{M}$ satisfies the $\matr{D}$-RIP with a constant $\delta_{k}$ then
\begin{eqnarray}
\label{eq:MP_RIP_norm}
\norm{\matr{M}\P_{T}}^2 \le 1+\delta_{k}
\end{eqnarray}
for every $T$ such that $\abs{T} \le k$.
\end{cor}

\begin{lem}
\label{lem:D_RIP_norm}
If $\;\matr{M}$ satisfies the $\matr{D}$-RIP~\eqref{def:D_RIP} then
\begin{eqnarray}
\label{eq:D_RIP_norm}
\norm{\P_{T}(\matr{I} - \matr{M}^*\matr{M})\P_{T}} \le \delta_k
\end{eqnarray}
for any $T$ such that $\abs{T} \le k$.
\end{lem}
{\em Proof:}
The proof is similar to the one of the standard RIP as appears in \cite{foucart10Sparse}.
We first observe that the definition~\eqref{def:D_RIP} of the $\matr{D}$-RIP is equivalent to requiring
$$
\abs{\norm{\matr{M}\vect{v}}_2^2  - \norm{\vect{v}}_2^2} \le \delta_{k}\norm{\vect{v}}_2^2
$$
for any $\vect{v} = \matr{D}\tilde{\alphabf}$ such that $\norm{\tilde\alphabf}_0\le k$. From this it follows that
$$
\abs{\norm{\matr{M}\P_{T}{\vect{z}}}_2^2  - \norm{\P_{T}{\vect{z}}}_2^2} \le \delta_{k}\norm{\P_{T}{\vect{z}}}_2^2 \le \delta_{k}\norm{{\vect{z}}}_2^2
$$
for any set $T$ such that $\abs{T} \le k$ and any ${\vect{z}} \in \Real^d$.
Next we notice that
\begin{eqnarray*}
\norm{\matr{M}\P_{T}{\vect{z}}}_2^2  - \norm{\P_{T}{\vect{z}}}_2^2 &=& {\vect{z}}^*\P_{T}\matr{M}^*\matr{M}\P_{T}{\vect{z}} - {\vect{z}}^*\P_{T}{\vect{z}}\\ \nonumber  &=& {\vect{z}}^*\P_{T}(\matr{M}^*\matr{M} - \matr{I}_d)\P_{T}{\vect{z}} \\ \nonumber
&=& \langle\P_{T}(\matr{M}^*\matr{M} - \matr{I}_d)\P_{T}{\vect{z}}, {\vect{z}} \rangle.
\end{eqnarray*}
Since $\P_{T}(\matr{M}^*\matr{M} - \matr{I}_d)\P_{T}$ is Hermitian we have that
\begin{eqnarray*}
\max_{{\vect{z}}}\frac{\langle\P_{T}(\matr{M}^*\matr{M} - \matr{I}_d)\P_{T}{\vect{z}}, {\vect{z}} \rangle}{\norm{{\vect{z}}}_2}
= \norm{\P_{T}(\matr{M}^*\matr{M} - \matr{I}_d)\P_{T}}.
\end{eqnarray*}
Thus we have that the $\matr{D}$-RIP implies \eqref{eq:D_RIP_norm} for any set $T$ such that $\abs{T}\le k$.
\hfill $\Box$ \bigskip

\begin{cor}
\label{cor:D_RIP_norm_diff}
If $\;\matr{M}$ satisfies the $\matr{D}$-RIP~\eqref{def:D_RIP} then
\begin{eqnarray}
\label{eq:D_RIP_norm_diff}
&& \norm{\P_{T_1}(\matr{I} - \matr{M}^*\matr{M})\P_{T_2}} \le \delta_k, \\ \nonumber
\end{eqnarray}
for any $T_1$ and $T_2$ with $\abs{T_1} \le k_1$, $\abs{T_2}\le k_2$, and $k_1+k_2\le k$.
\end{cor}
{\em Proof:}
Since $T_1 \subset T_1 \cup T_2$ and $T_2 \subset T_1 \cup T_2$, we have from Lemma~\ref{lem:D_RIP_norm} that
\begin{eqnarray}
\nonumber
\norm{\P_{T_1}(\matr{I} - \matr{M}^*\matr{M})\P_{T_2}}_2 \le \norm{\P_{T_2 \cup T_1}(\matr{I} - \matr{M}^*\matr{M})\P_{T_2 \cup T_1}} \leq \delta_k.
\end{eqnarray}

\hfill $\Box$ \bigskip

Finally, we point out some consequences of the definition of near-optimal projections, as in Definition~\ref{def:C_optimal_proj}.  A clear implication of this definition is that for any vector $\vect{v}\in \Real^d$ that has a $k$-sparse representation and a support set $T$ such that $\abs{T} \le k$,
and for any $\vect{z} \in \RR{d}$  we have that
\begin{eqnarray}
\label{eq:C_optimal_ineq}
&& \norm{\vect{z}-\P_{\CFI_{\pscI k}(\vect{z})}\vect{z}}_2^2 \le \CI_k \norm{\vect{v}-\vect{z}}_2^2, \quad\text{and}\\
\label{eq:C_optimal_ineq_up}
&& \norm{\P_{\CFII_{\pscII k}(\vect{z})}\vect{z}}_2^2 \ge {\CII}_k\norm{\P_{T}\vect{z}}_2^2.
\end{eqnarray}

The constants $\CI_k$ and ${\CII}_{2k}$ will play a role in the convergence guarantees we develop for SSCoSaMP.
Requirements on the allowed values and the type of dictionaries that has near optimal support selection schemes will be discussed later in Section~\ref{sec:near_opt_proj}.
We will also utilize the following elementary fact, whose proof is immediate using the inequality of arithmetic and geometric means.
\begin{prop}
\label{prop:norm2_ineq}
For any two given vectors $\vect{x}_1$, $\vect{x}_2$ and a constant $c>0$ it holds that
\begin{eqnarray}
\norm{\vect{x}_1+\vect{x}_2}_2^2 \le (1+c)\norm{\vect{x}_1}_2^2 + \left(1+\frac{1}{c}\right)\norm{\vect{x}_2}_2^2.
\end{eqnarray}
\end{prop}


\section{Algorithm Guarantees}
\label{sec:guarantees}

In this section we provide theoretical guarantees for the reconstruction performance of SSCoSaMP. 
The results here are for the choice of $a = 2$ in the algorithm, however, analogous results for other values $a\geq 1$ follow similarly. 
We will prove the main result, Theorem~\ref{thm:general_bound}, via Corollary~\ref{cor:SSCoSaMP_bound}.  The proof and discussion of this corollary occupy the remainder of this section.

\subsection{Theorem Conditions}

Before we begin the proof of the theorem we first ask under what conditions the assumptions of the theorem hold.
One condition of Theorem~\ref{thm:general_bound} is that $\delta_{2(1+\pscI)k} \le \epsilon^2_{\CI_k,\CII_{2k},\gamma}$
for a constant $\epsilon^2_{\CI_k,\CII_{2k},\gamma} >0$.
When the dictionary $\matr{D}$ is unitary, it was shown for many families of random matrices that for any value of $\epsilon_k$,
if $m \ge \frac{C}{\epsilon_k^2} k \log(\frac{m}{k\epsilon_k})$, where $C$ is a given constant, then $\delta_k \le \epsilon_k$ with high probability~\cite{Candes06Near, Rauhut08Compressed,Mendelson08Uniform}.
A similar result for the same family of random matrices holds for the $\matr{D}$-RIP \cite{Candes11Compressed}.
Thus, the critical part in the conditions of the Theorem is condition \eqref{eq:C_k_tilda_C_2k_cond}, that imposes a requirement on $\CI_k$ and $\CII_{2k}$ to be close to $1$.
We have an access to projection operators that satisfy this condition in many practical settings
which are not supported by the guarantees provided in previous papers that used near optimal projections \cite{Davenport13Signal,Blumensath11Sampling,Giryes13Greedy}.
This is due to the near-optimality definition and the proof technique used in this paper; A detailed discussion of this subject is left to Section~\ref{sec:near_opt_proj} below.

\subsection{SSCoSaMP Guarantees}

Analogously to that of CoSaMP in \cite{foucart10Sparse}, our proof relies on iteration invariant which shows that each iteration substantially reduces the recovery error.

\begin{thm}
\label{thm:SSCoSaMP_iter_bound}
Let $\vect{y} = \matr{M}\x +\vect{e}$, where $\matr{M}$ satisfies the $\matr{D}$-RIP~\eqref{def:D_RIP} with a constant $\delta_{(3\pscI+1)k}$ ($\pscI\ge 1$), $\x$ is a vector with a $k$-sparse representation under $\matr{D}$ and $\vect{e}$ is an additive noise vector.
Suppose $\CFI_{\pscI k}$ and $\CFII_{2\pscI k}$ be near optimal projections as in Definition~\ref{def:C_optimal_proj} with constants $\CI_k$ and ${\CII}_{2k}$.
Then
\begin{eqnarray}
\label{eq:SSCoSaMP_iter_bound}
&& \hspace{-0.3in} \norm{\vect{x}^t -\vect{x}}_2 \le \rho\norm{\vect{x} - \vect{x}^{t-1}}_2 +
  \eta\norm{\vect{e}}_2,
\end{eqnarray}
for constants $\rho$ and $\eta$.
The iterates converge, i.e. $\rho<1$, if $\delta_{(3\pscI+1)k}< \epsilon^2_{\CI_k,{\CII}_{2k},\gamma}$, for some positive constant $\epsilon^2_{\CI_k,{\CII}_{2k},\gamma}$, and \eqref{eq:C_k_tilda_C_2k_cond} holds.
\end{thm}

An immediate corollary of the above theorem is the following
\begin{cor}
\label{cor:SSCoSaMP_bound}
Assume the conditions of Theorem~\ref{thm:SSCoSaMP_iter_bound}.  Then after a constant number of iterations $t^* = \ceil{\frac{\log(\norm{\vect{x}}_2/\norm{\vect{e}}_2)}{\log(1/\rho)}}$ it holds that
\begin{eqnarray}
\label{eq:SSCoSaMP_bound}
&& \hspace{-0.5in} \norm{\vect{x}^{t^*} -\vect{x}}_2 \le
  \left(1 +  \frac{1-\rho^{t^*}}{1-\rho}\right)\eta\norm{\vect{e}}_2.
\end{eqnarray}
\end{cor}
{\em Proof:}
By using \eqref{eq:SSCoSaMP_iter_bound} and recursion we have that after $t^*$ iterations
\begin{eqnarray}
&& \hspace{-0.3in} \norm{\vect{x}^{t^*} -\vect{x}}_2 \le \rho^{t^*}\norm{\vect{x} - \vect{x}^{0}}_2 +
  (1+\rho+\rho^2+\dots \rho^{t^*-1})\eta\norm{\vect{e}}_2
  \\ \nonumber && \le\left(1 +  \frac{1-\rho^{t^*}}{1-\rho}\right)\eta\norm{\vect{e}}_2,
\end{eqnarray}
where the last inequality is due to the equation of the geometric series, the choice of $t^*$, and the fact that $\vect{x}^0 =\vect{0}$.
\hfill $\Box$ \bigskip

Note that Corollary~\ref{cor:SSCoSaMP_bound} implies our main result, Theorem~\ref{thm:general_bound}, with $\eta_0 = \left(1 +  \frac{1-\rho^{t^*}}{1-\rho}\right)\eta$.

We turn now to prove the iteration invariant, Theorem~\ref{thm:SSCoSaMP_iter_bound}.  Instead of presenting the proof directly, we divide the proof into several lemmas.
The first lemma gives a bound for $\norm{\vect{x}_p -\vect{x} }_2$ as a function of $\norm{\vect{e}}_2$ and $\norm{\Q_{\tilde{T}^t}(\vect{x}_p - \vect{x})}_2$.
\begin{lem}
\label{lem:SSCoSaMP_xp_bound}
If $\matr{M}$ has the $\matr{D}$-RIP with a constant $\delta_{3\pscI k}$, then with the notation of Algorithm~\ref{alg:SSCoSaMP}, we have
\begin{eqnarray}
\label{eq:SSCoSaMP_xp_bound}
\norm{\vect{x}_p -\vect{x}}_2 \le \frac{1}{\sqrt{1-\delta_{(3\pscI +1)k}^2}}\norm{\Q_{\tilde{T}^t}(\vect{x}_p - \vect{x})}_2 
+  \frac{\sqrt{1+\delta_{3\pscI k}}}{1-\delta_{(3\pscI +1)k}}\norm{\vect{e}}_2
\end{eqnarray}
\end{lem}
The second lemma bounds $\norm{\vect{x}^{t} - \vect{x}}_2$ in terms of $\norm{\Q_{\tilde{T}^t}(\vect{x}_p - \vect{x})}_2$ and $\norm{\vect{e}}_2$ using the first lemma.
\begin{lem}
\label{lem:SSCoSaMP_xt_bound1}
Under the assumptions and notation of Theorem~\ref{thm:general_bound}, we have
\begin{eqnarray}
&& \hspace{-0.5in} \norm{\vect{x}^t -\vect{x}}_2 \le \rho_1\norm{\Q_{\tilde{T}^t}(\vect{x}_p - \vect{x})}_2+  \eta_1\norm{\vect{e}}_2,
\end{eqnarray}
where the constants $\rho_1$ and $\eta_1$ are given explicitly in~\eqref{constants}.
\end{lem}
The last lemma bounds $\norm{\Q_{\tilde{T}^t}(\vect{x}_p - \vect{x})}_2$ with $\norm{\vect{x}^{t-1} - \vect{x}}_2$ and $\norm{\vect{e}}_2$.
\begin{lem}
\label{lem:SSCoSaMP_Pxp_bound}
Under the assumptions and notation of Theorem~\ref{thm:general_bound}, we have
\begin{eqnarray}
\label{eq:SSCoSaMP_Pxp_bound}
&& \hspace{-0.3in}\norm{\Q_{\tilde{T}^t}(\vect{x}_p - \vect{x})}_2 \le 
 \eta_2\norm{\vect{e}}_2
 + \rho_2\norm{\vect{x} - \vect{x}^{t-1}}_2,
\end{eqnarray}
where the constants $\rho_2$ and $\eta_2$ are given explicitly in~\eqref{constants}.
\end{lem}
The proofs of Lemmas~\ref{lem:SSCoSaMP_xp_bound}, \ref{lem:SSCoSaMP_xt_bound1} and \ref{lem:SSCoSaMP_Pxp_bound} appear in \ref{sec:SSCoSaMP_xp_bound_proof}, \ref{sec:SSCoSaMP_xt_bound1_proof} and \ref{sec:SSCoSaMP_Pxp_bound_proof}, respectively. 
With the aid of the above three lemmas we turn to the proof of the iteration invariant, Theorem~\ref{thm:SSCoSaMP_iter_bound}.

{\em Proof of Theorem~\ref{thm:SSCoSaMP_iter_bound}:}
Substituting the inequality of Lemma~\ref{lem:SSCoSaMP_Pxp_bound} into the inequality of Lemma~\ref{lem:SSCoSaMP_xt_bound1} gives
\eqref{eq:SSCoSaMP_iter_bound} with $\rho = \rho_1\rho_2$ and $\eta = \eta_1 + \rho_1\eta_2$. The iterates converge if $\rho_1^2\rho_2^2 < 1$.
Since $\delta_{(\pscI+1)k} \le \delta_{3\pscI k} \le \delta_{(3\pscI+1)k}$ this holds if
\begin{eqnarray}\label{interinequ}
&& \hspace{-0.3in} \frac{\left(1+\sqrt{\CI_k} \right)^2 }{1-\delta_{(3\pscI+1)k}^2} 
   \left(1-\left( \left(\frac{\sqrt{\CII_{2k}}}{1+\gamma} +1 \right)\sqrt{\delta_{(3\pscI+1)k}} -\frac{\sqrt{\CII_{2k}}}{1+\gamma}
 \right)^2\right) < 1.
\end{eqnarray}
Since $\delta_{(3\pscI+1)k} <1$, we have $\delta_{(3\pscI+1)k}^2 < \delta_{(3\pscI+1)k} $.  Using this fact and expanding~\eqref{interinequ} yields
the stricter condition
\begin{eqnarray}
\label{eq:epsilon_quadratic_ineq}
&&  \hspace{-0.3in} 
\left( \frac{1}{\left(1+\sqrt{\CI_k} \right)^2} -
 \left(\frac{\sqrt{\CII_{2k}}}{1+\gamma} +1 \right)^2\right) {\delta_{(3\pscI+1)k}} +
 2\left(\frac{\sqrt{\CII_{2k}}}{1+\gamma} +1 \right)\frac{\sqrt{\CII_{2k}}}{1+\gamma}\sqrt{\delta_{(3\pscI+1)k}} 
\\ \nonumber && \hspace{3in} +
1 -  \frac{1}{\left(1+\sqrt{\CI_k} \right)^2} -\frac{{\CII_{2k}}}{\left(1+\gamma\right)^2}  < 0.
\end{eqnarray}
The above equation has a positive solution if and only if \eqref{eq:C_k_tilda_C_2k_cond} holds.
Denoting its positive solution by ${\epsilon_{\CI_k,\CII_{2k},\gamma}}$ we have that the expression holds when
$\delta_{(3\pscI+1)k} \le \epsilon_{\CI_k,\CII_{2k},\gamma}^2$, which completes the proof.  Note that in the proof we have 
\begin{eqnarray}\label{constants}
 {\eta_1 = \frac{\left(1 + \sqrt{C_k}\right)\sqrt{1+\delta_{3\pscI k}}}{1-\delta_{(3\psc +1)k}}}, &&
{\eta_2^2 = \bigg(\frac{1+\delta_{3\pscI k}}{\gamma(1+\alpha)}
 +\frac{(1+\delta_{(\pscI +1)k})\CII_{2k}}{\gamma(1+\alpha)(1+\gamma)}\bigg)},\\
  {\rho_1^2 = \frac{\left( 1+\sqrt{\CI_k} \right)^2}{1-\delta_{(3\pscI+1)k}^2}}, &&
    {\rho_2^2 = 1-\bigg( \sqrt{\delta_{(3\pscI+1)k}}-\frac{\sqrt{\CII_{2k}}}{1+\gamma}\left(1-\sqrt{\delta_{(\pscI+1)k}}\right)
 \bigg)^2 }, \notag\\ &&
   {\alpha = \frac{\sqrt{\delta_{(3\pscI+1)k}}}{\sqrt{\frac{\CII_{2k}}{(1+\gamma_1)(1+\gamma_2)}}\left(1-\sqrt{\delta_{(\pscI+1)k}}\right)- \sqrt{\delta_{(3\pscI+1)k}}}}\notag
   \end{eqnarray}
     and $\gamma >0$ is an arbitrary constant.
\hfill $\Box$ \bigskip


\section{Near Optimal Projection Examples}
\label{sec:near_opt_proj}

In this section we give several examples for which condition \eqref{eq:C_k_tilda_C_2k_cond}, 
 \begin{eqnarray*}
\left(1+\sqrt{\CI_k} \right)^2\left( 1-\frac{{\CII}_{2k}}{(1+\gamma)^2} \right)<1,
\end{eqnarray*}
can be satisfied with accessible projection methods.

\subsection{Unitary Dictionaries}
For unitary $\matr{D}$ the conditions hold trivially since $\CI_k = {\CII}_{2k} = 1$
using simple thresholding.
In this case our results coincide with the standard representation model for which we already have theoretical guarantees~\cite{Needell09CoSaMP}.
However, for a general dictionary $\matr{D}$ simple thresholding is not expected to have this property.

\subsection{RIP Dictionaries}
We next consider the setting in which the dictionary $\matr{D}$ itself satisfies the RIP.
In this case we may use a standard method like IHT or CoSaMP for $\CFI_{k}$ and simple thresholding for $\CFII_{2k}$. 
For dictionaries that satisfy the RIP, it is easy to use existing results in order to derive bounds on the constant $\CI_k$.

In order to see how such bounds can be achieved, notice that standard bounds exists for these techniques in terms of the representation error rather than the signal error.
That is, for a given vector $\vect{v} = \matr{D}\alphabf + \vect{e}$ and any support set $T^*$ of size $k$,
it is guaranteed that if $\delta_{4k} \le 0.1$ (or $\delta_{3k} \le \frac{1}{\sqrt{32}}$) then $\hat{\alphabf}$, the recovered representation of CoSaMP (or IHT),
  satisfies
\begin{eqnarray}
\label{eq:greedy_repres_recovery_bound}
\norm{\alphabf - \hat{\alphabf}}_2 \le C_{e} \norm{\vect{v} - \matr{P}_{T^*}\vect{v}}_2,
\end{eqnarray}
where $C_e \simeq 5.6686$ (or $C_e \simeq 3.3562$) \cite{foucart10Sparse, Needell09CoSaMP, Blumensath09Iterative}.

We use this result to bound $\CI_k$ as follows.
For a general vector $\vect{v}$ we may write its optimal projection as $\matr{P}_{T^*}\vect{v} = \matr{D}\alphabf$ with $\supp(\vect{\alpha}) = T^*$.  Applying the bound in \eqref{eq:greedy_repres_recovery_bound} with $\vect{e} = \vect{v} - \matr{P}_{T^*}\vect{v}$ and $\matr{P}_{\hat{T}}\vect{v} = \matr{D}\vect{\hat{\alphabf}}$ along with the RIP yields

\begin{align}
\norm{\vect{v} - \matr{P}_{\hat{T}}\vect{v}}_2 &\leq \norm{\vect{v} - \matr{P}_{T^*}\vect{v}}_2 + \norm{\matr{P}_{T^*}\vect{v} - \matr{P}_{\hat{T}}\vect{v}}_2\\
\nonumber
 &= \norm{\vect{v} - \matr{P}_{T^*}\vect{v}}_2 + \norm{\matr{D}\alphabf - \matr{D}\vect{\hat{\alphabf}}}_2\\
\nonumber
&\leq \norm{\vect{v} - \matr{P}_{T^*}\vect{v}}_2 + \sqrt{1 + \delta_{2k}}\norm{\alphabf - \vect{\hat{\alphabf}}}_2\\
\nonumber
&\leq \norm{\vect{v} - \matr{P}_{T^*}\vect{v}}_2 + C_e\sqrt{1+\delta_{2k}}\norm{\vect{v} - \matr{P}_{T^*}\vect{v}}_2.
\end{align}

This implies that
\begin{eqnarray}
\label{eq:C_k_C_e_bound}
\CI_k \le 1+C_e\sqrt{1+\delta_{2k}}.
\end{eqnarray}
For example, if $\delta_{4k} \le 0.1$ then $\CI_k \le 6.9453$ for CoSaMP and if $\delta_{3k}\le \frac{1}{\sqrt{32}}$ then $\CI_k \le 4.6408$ for IHT.
The inequality in \eqref{eq:C_k_C_e_bound} holds true not only for CoSaMP and IHT but for any algorithm that provides a $k$-sparse representation that obeys the bound in \eqref{eq:greedy_repres_recovery_bound}.
Note that many greedy algorithms have these properties (e.g. \cite{Foucart11Hard, Dai09Subspace, Zhang11Sparse}), but relaxation techniques such as $\ell_1$-minimization \cite{Candes05Decoding} or the Dantzig selector \cite{Candes07Dantzig} are not guaranteed to give a $k$-sparse result.

Having a bound for $\CI_k$, we realize that in order to satisfy \eqref{eq:C_k_tilda_C_2k_cond} we now have a condition on the second constant, 
\begin{eqnarray}
{\CII}_{2k} \ge \left( 1- \frac{1}{\left(1+\sqrt{\CI_k}\right)^2}\right)(1+\gamma)^2.
\end{eqnarray}
In order to show that this condition can be satisfied we provide an upper bound for ${\CII}_{2k}$
which is a function of the RIP constants of $\matr{D}$.  The near-optimal projection can be obtained by simple thresholding under the image of $\matr{D}^*$:
\begin{equation}\label{Dthres}
\CFII_{k}(\vect{v}) = \argmin_{|T|=k} \|\matr{D}^*_T\vect{v}\|_2.
\end{equation}

\begin{lem}[Thresholding Projection RIP bound]
If $\;\matr{D}$ is a dictionary that satisfies the RIP with a constant $\delta_k$, then
using~\eqref{Dthres} as the thresholding projector yields $${\CII}_{k} \ge \frac{1-\delta_{k}}{1+\delta_{k}}.$$
\end{lem}
{\em Proof:}
Let $\vect{v}$ be a general vector. Let $\hat{T}$ be the indices of the largest $k$ entries of $\matr{D}^*\vect{v}$
and $T^*$ the support selected by the optimal support selection scheme as in~\eqref{eq:optimal_sparse_projection}.
By definition we have that
\begin{eqnarray}
\norm{\matr{D}_{\hat{T}}^*\vect{v}}_2^2 \ge \norm{\matr{D}_{T^*}^*\vect{v}}_2^2.
\end{eqnarray}
Since $\frac{1}{1+\delta_k}\le \norm{(\matr{D}_T^*)^\dag}_2^2 \le \frac{1}{1-\delta_k}$
for $\abs{T}\le k$ (see Prop. 3.1 of~\cite{Needell09CoSaMP}),
we have that
\begin{eqnarray}
\left(1+\delta_k\right)\norm{(\matr{D}_{\hat{T}}^*)^\dag\matr{D}_{\hat{T}}^*\vect{v}}_2^2 \ge \left(1-\delta_k\right)\norm{(\matr{D}_{T^*}^*)^\dag\matr{D}_{T^*}^*\vect{v}}_2^2.
\end{eqnarray}
Since $\P_{\hat{T}} = (\matr{D}_{\hat{T}}^*)^\dag\matr{D}_{\hat{T}}^*$ we get that
\begin{eqnarray}
\norm{\P_{\hat{T}}\vect{v}}_2^2 \ge \frac{1-\delta_k}{1+\delta_k}\norm{\P_{T^*}\vect{v}}_2^2.
\end{eqnarray}
Thus ${\CII}_k \ge \frac{1-\delta_k(\matr{D})}{1+\delta_k(\matr{D})}$.
\hfill $\Box$ \bigskip

Hence, the condition on the RIP of $\matr{D}$ for satisfying \eqref{eq:C_k_tilda_C_2k_cond} turns to be 
\begin{eqnarray}\label{constant2}
\frac{1-\delta_{2k}}{1+\delta_{2k}} \ge \left( 1- \frac{1}{\left(1+\sqrt{\CI_k}\right)^2}\right)(1+\gamma)^2.
\end{eqnarray}

By using the exact expression for $\CI_k$ in terms of RIP constants, one can obtain guarantees in terms of the RIP constants only.  For example, from~\cite{Needell09CoSaMP}, for CoSaMP one has more precisely that for any vector $\alphabf$, the reconstructed vector $\hat{\alphabf}$ from measurements $\vect{z} = \matr{D}{\alphabf} + \vect{e}$ satisfies
$$
\norm{\vect{\alphabf} - \vect{\hat{\alphabf}}} \leq \left[\frac{2}{\sqrt{1-\delta_{3k}}} + 4\left(1 + \frac{\delta_{4k}}{1-\delta_{3k}}\right)\cdot\frac{1}{\sqrt{1-\delta_{2k}}}\right]\norm{\vect{e}}.
$$

Using \eqref{eq:C_k_C_e_bound} we have
\begin{eqnarray}
\CI_k \leq 1 + \sqrt{1+\delta_{2k}}\left[\frac{2}{\sqrt{1-\delta_{3k}}} + 4\left(1 + \frac{\delta_{4k}}{1-\delta_{3k}}\right)\cdot\frac{1}{\sqrt{1-\delta_{2k}}}\right].
\end{eqnarray}
Substituting this into the expression~\eqref{constant2} gives a bound on the RIP constants alone.  For example, setting $\gamma = 0.01$, one finds that the requirement $\delta_{4k} \leq 0.027$ is enough to guarantee~\eqref{eq:C_k_tilda_C_2k_cond} holds using CoSaMP.

\subsection{Incoherent Dictionaries}

Given that a  dictionary $\matr{D}$ has a coherence $\mu$,
it is known that the RIP constant can be upper bounded by $\mu$ in the following way \cite{Elad10Sparse}
\begin{eqnarray}
\delta_k \le (k-1)\mu.
\end{eqnarray}
Hence, using this relation one may get recovery conditions based on the coherence value using the conditions from the previous subsection.
For example, if we use CoSaMP for the first projection and thresholding for the second one, one may have the following condition in terms of the coherence (instead of the RIP): 
$\mu \leq \frac{0.027}{4k-1}$.

\subsection{Support Selection using Highly Correlated Dictionaries}
\label{sec:high_cor}
In all the above cases, the dictionary is required to be incoherent.
This follows from the simple fact that decoding under a coherent dictionary is a hard problem in general. However, in some cases we have a coherent dictionary in which each atom has a high correlation with a small number of other atoms and very small correlation with all the rest.
In this case, the high coherence is due to these rare high correlations
and pursuit algorithms may fail to select the right atoms in their support estimate as they may be confused between the right atom and its highly correlated columns.
Hence, one may update the pursuit strategies to add in each of their steps only
atoms which are not highly correlated with the current selected atoms and as a final stage extend the estimated support to include all the atoms which have high coherence with the selected support set.

This idea is related to the recent literature of super-resolution (see e.g.~\cite{schmidt1986MUSIC,fannjiang2012CS,candes2013towards,
demanet2013super,divekar2013super,Duarte13Spectral} and references therein) 
and to the $\epsilon$-OMP algorithm \cite{Giryes13OMP}, which is an extension of OMP.
In this work we employ $\epsilon$-OMP (with a post-processing step that adds correlated atoms) as a support selection procedure.
We also propose a similar extension for thresholding, $\epsilon$-thresholding, that for a given signal $\vect{z}$, selects the support in the following way. It picks the indices of the largest elements of $\matr{D}^*\vect{z}$ one at a time, where at each time it adds the atom with highest correlation to $\vect{z}$ excluding the already selected ones. Each atom is added together with its highly correlated columns.

Before we present these methods formally, we introduce the following definition taken from \cite{Giryes13OMP}.
\begin{defn}[$\epsilon$-extension\footnote{In \cite{Giryes13OMP} it is referred to as $\epsilon$-closure but since closure bears a different meaning in mathematics we use a different name here.}]
Let $0\le \epsilon < 1$ and $\matr{D}$ be a fixed dictionary.
The $\epsilon$-extension of a given support set $T$ is defined as
$$
\ext_{\epsilon,2}(T) = \left\{i\; : \;\exists j\in T, ~ \frac{\abs{\langle\matr{d}_i, \matr{d}_j\rangle}^2}{\norm{\matr{d}_i}_2^2\norm{\matr{d}_j}_2^2}\ge  1-\epsilon^2\right\}.
$$
\end{defn}
Having the above definition, we present $\epsilon$-OMP\footnote{In \cite{Giryes13OMP} $\epsilon$-OMP is presented slightly different: (1) It treats the more general case of recovering a signal from a set of measurement $\vect{y} = \matr{M}\vect{x} + \vect{e}$; (2)
the support extension at the last stage of Algorithm~\ref{alg:OMP_eps} is proposed as a post-processing step apart from the $\epsilon$-OMP algorithm.} and $\epsilon$-Thresholding techniques in Algorithms~\ref{alg:OMP_eps}  and \ref{alg:thresh_eps}.

 \begin{algorithm}[t]
 \caption{$\epsilon$-Orthogonal Matching Pursuit} \label{alg:OMP_eps}
\begin{algorithmic}[l]

\REQUIRE $k, \matr{D}, \vect{z}$ where $\vect{z} = \vect{x}
+ \vect{e}$, $\vect{x} = \matr{D}\alphabf$, $\norm{\alphabf}_0 \le k$
and $\vect{e}$ is additive noise.

\ENSURE $\hat{\vect{x}}$: $k$-sparse approximation of
$\vect{x}$ supported on $\hat{T}$.

\STATE Initialize estimate $\hat{\x}^0 = \vect{0}$, residual $\vect{r}^0 = \vect{z}$, support $\hat{T}^0 =\check{T}^0 =\emptyset$
 and set $t = 0$.

\WHILE{$t \le k$}

\STATE $t = t + 1$.

\STATE New support element: $i^t = \argmax_{i \not \in \check{T}^{t-1}}
|\matr{d}^*_i\vect{r}^{t-1}|$.

\STATE Extend support: $\hat{T}^t = \hat{T}^{t-1} \cup \{i^t\}$.

\STATE Calculate a new estimate: $\hat{\vect{x}}^t = \matr{D}_{\hat{T}^t}\matr{D}_{\hat{T}^t}^\dag\vect{z}$.

\STATE Calculate a new residual: $\vect{r}^t = \vect{z} - \hat{\vect{x}}^t$.

\STATE Support $\epsilon$-extension: $\check{T}^t = \ext_{\epsilon,2}(\hat{T}^{t})$.

\ENDWHILE

\STATE Set estimated support $\hat{T} = \check{T}^t$.

\STATE Form the final solution $\hat{\vect{x}} =  \matr{D}_{\hat{T}}\matr{D}_{\hat{T}}^\dag\vect{z}$.

\end{algorithmic}
\end{algorithm}

 \begin{algorithm}[t]
\caption{$\epsilon$-thresholding} \label{alg:thresh_eps}
\begin{algorithmic}[l]

\REQUIRE $k, \matr{D}, \vect{z}$ where $\vect{z} = \vect{x}
+ \vect{e}$, $\vect{x} = \matr{D}\alphabf$, $\norm{\alphabf}_0 \le k$
and $\vect{e}$ is additive noise.

\ENSURE $\hat{\vect{x}}$: a $k$-sparse approximation of
$\vect{x}$ supported on $\hat{T}$.

\STATE Initialize support $\hat{T}^0 =\check{T}^0 =\emptyset$ and set $t = 0$.

\STATE Calculate correlation between dictionary and measurements: $\vect{v} =
\matr{D}^*\vect{z}$.

\WHILE{$t \le k$}

\STATE $t = t + 1$.

\STATE New support element: $i^t = \argmax_{i \not \in \check{T}^{t-1}}
|\vect{v}_i|$.

\STATE Extend support: $\hat{T}^t = \hat{T}^{t-1} \cup \{i^t\}$.

\STATE Support $\epsilon$-extension: $\check{T}^t = \ext_{\epsilon,2}(\hat{T}^{t})$.

\ENDWHILE

\STATE Set estimated support $\hat{T} = \check{T}^t$.

\STATE Form the final solution $\hat{\vect{x}} =  \matr{D}_{\hat{T}}\matr{D}_{\hat{T}}^\dag\vect{z}$.
\end{algorithmic}
\end{algorithm}

Note that the size of the group of atoms which are highly correlated with one atom of $\matr{D}$ is bounded. The size of the largest group is an upper bound for the near-optimality constants $\pscI$ and $\pscII$ (note here we will just set $\pscI = \pscII$). More precisely, if the allowed high correlations
are greater then $1-\epsilon^2$ then we have the upper bound $$\pscI \le
\max_{T:\abs{T} \le k} \abs{\ext_{\epsilon,2}(T)} \le \max_{1 \le i \le n} k\abs{\ext_{\epsilon,2}(\{i\})}.$$

We have a trade-off between the size of the correlation which we can allow and the size of the estimated support which we get. The smaller the correlation between columns we allow, the larger $\pscI$ is and thus also the estimated support.
On the one hand, this attribute is positive; the larger the support, the higher the probability that  our near-optimality constants $\CI_k$ and $\CII_k$ are close to $1$.
On the other hand, for large $\pscI$, $\delta_{(3\pscI+1) k}$ is larger and it is harder to satisfy the RIP requirements.
Hence we expect that if the number of measurements is small, the size of $\pscI$ would be more critical as it would be harder to satisfy the RIP condition. When the number of measurements gets higher, the RIP requirement is easier to satisfy and can handle higher values of $\pscI$.

One trivial example, in which the above projections have known near-optimality constants is when $\matr{D}$ is an incoherent dictionary with one or more repeated columns.
In this case, the projection constants of $\matr{D}$ are simply the ones of the underlying incoherent dictionary.

\begin{figure}[htb]
  \centering
  \centerline{\includegraphics[width=5.0cm]{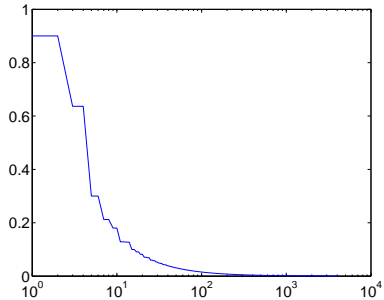}}
\caption{Correlation size (inner product) in a sorted order of one atom of the $4$ times redundant-DFT dictionary with the other atoms.
Note that the x-axis is in a log-scale.}
\label{fig:DFT4_correlations}
\end{figure}

In other cases we still do not have guarantees for these constants.
In \cite{Giryes13OMP}, a more general version of $\epsilon$-OMP that includes the matrix $\matr{M}$ is analyzed, providing conditions for the signal reconstruction. However, these impose requirements on the magnitude of the signal coefficients which we do not have control of in the projection problem.
Hence, the existing recovery guarantees  for $\epsilon$-OMP cannot be used
for developing bounds for the projection constants.

Though theoretical statements are not at hand yet, we shall see that
these methods give good recovery in practice.
Clearly, we need each atom in $\matr{D}$ to be highly correlated only with
a small group of other columns and incoherent with all the rest.
An example of such a dictionary is the overcomplete-DFT which is a highly coherent dictionary.
The correlations between each atom in this dictionary and its neighboring atoms are the same, i.e., each of the diagonals of its Gram matrix have the same value. A plot of the coherence value of a given atom with its neighbors in a sorted order appears in Fig.~\ref{fig:DFT4_correlations} for a four times overcomplete DFT and signal dimension $d=1024$.

Note that when we determine a correlation to be high, if the inner product (atoms are normalized) between two atoms is greater than $0.9$ ($\epsilon = \sqrt{0.1}$), we get that each atom has two other highly correlated columns with correlation of size $0.9$. The correlation with the rest is below $0.64$, where the largest portion has inner products smaller then $0.1$.

\subsubsection{Experimental Results}

We repeat the experiments from \cite{Davenport13Signal} for the overcomplete-DFT with redundancy factor $4$ and check the effect of the new support selection methods
both for the case where the signal coefficients are clustered and the case where they are well separated.  We compare the performance of OMP, $\epsilon$-OMP and $\epsilon$-thresholding for the approximate projections. We do not include other methods since a thorough comparison has been already performed in \cite{Davenport13Signal}, and the goal here is to check the effect of the $\epsilon$-extension step.

\begin{figure}[htb]
\begin{minipage}[b]{.48\linewidth}
  \centering
  \centerline{\includegraphics[width=7.0cm]{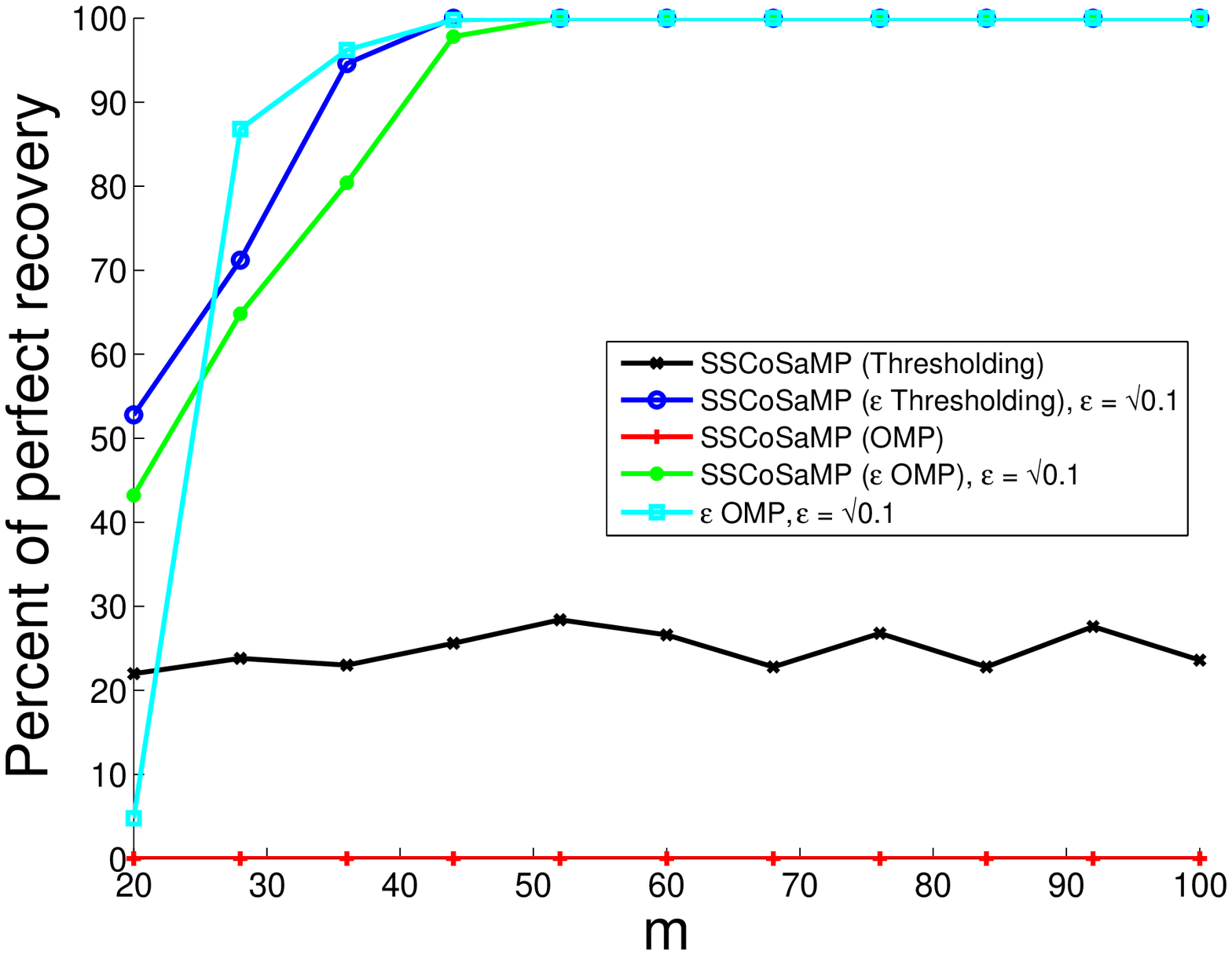}}
\end{minipage}
\hfill
\begin{minipage}[b]{.48\linewidth}
  \centering
  \centerline{\includegraphics[width=7.0cm]{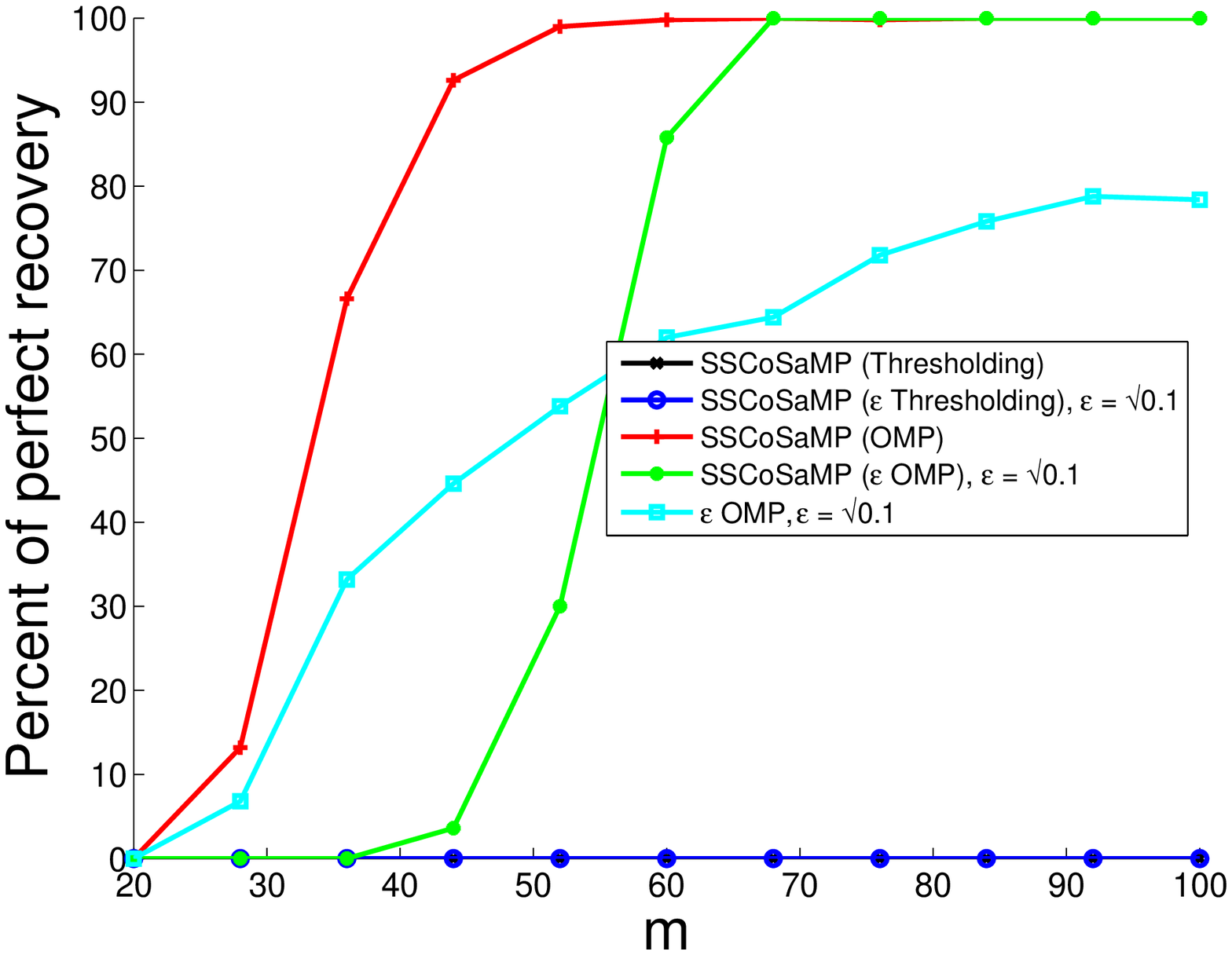}}
\end{minipage}
\caption{Recovery rate for SSCoSaMP (Thresholding), SSCoSaMP ($\epsilon$-Thresholding) with $\epsilon = \sqrt{0.1}$,
SSCoSaMP (OMP), SSCoSaMP ($\epsilon$-OMP) with $\epsilon = \sqrt{0.1}$
and $\epsilon$-OMP with $\epsilon = \sqrt{0.1}$ for a random $m\times 1024$ Gaussian matrix $\matr{M}$
and a $4$ times overcomplete DFT matrix $\matr{D}$.  The signal is $8$-sparse and on the left the coefficients of the original signal are clustered whereas on the right they are separated. 
}
\label{fig:recovery rate}
\end{figure}

\begin{figure}[htb]
\begin{minipage}[b]{.48\linewidth}
  \centering
  \centerline{\includegraphics[width=7.0cm]{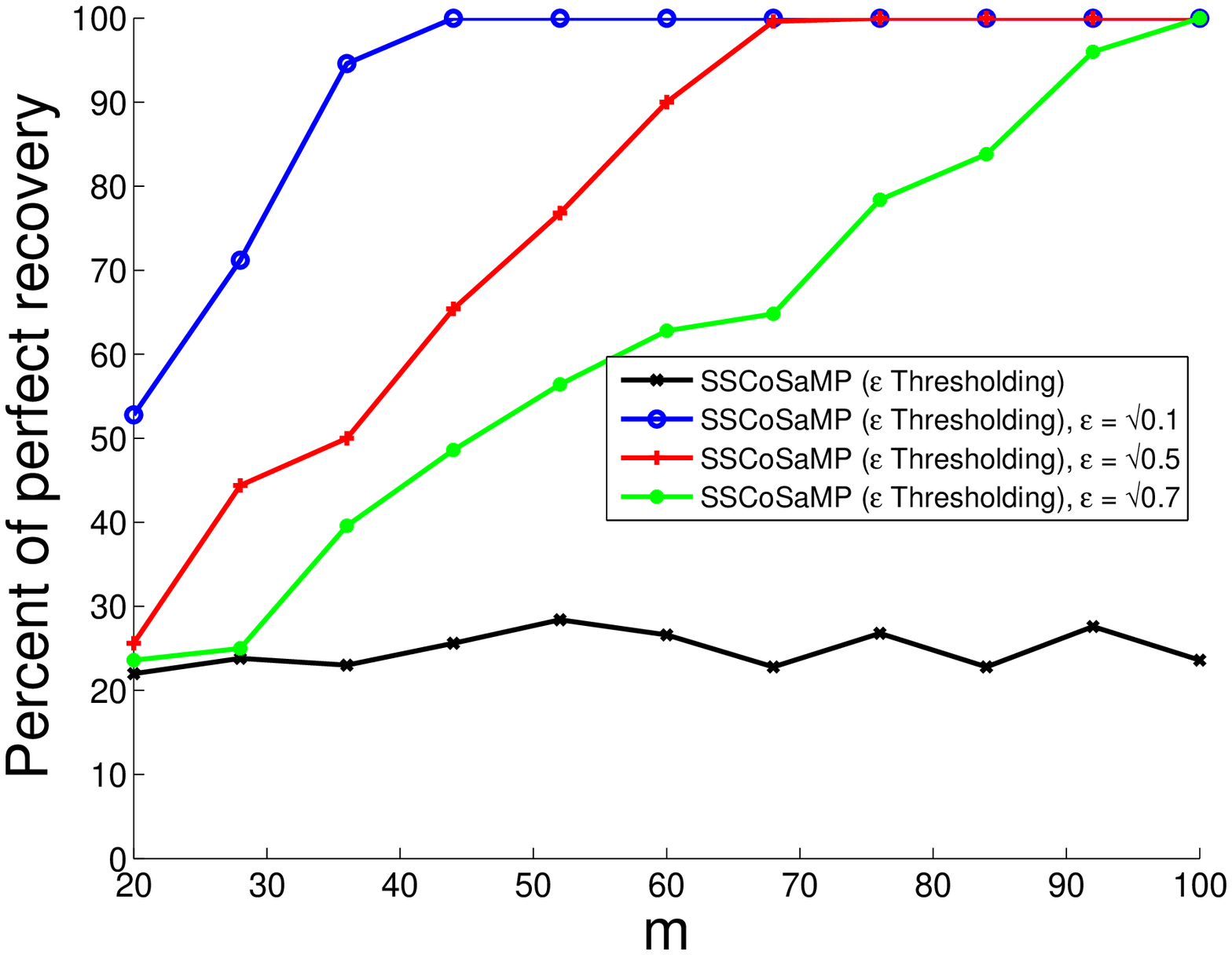}}
\end{minipage}
\hfill
\begin{minipage}[b]{.48\linewidth}
  \centering
  \centerline{\includegraphics[width=7.0cm]{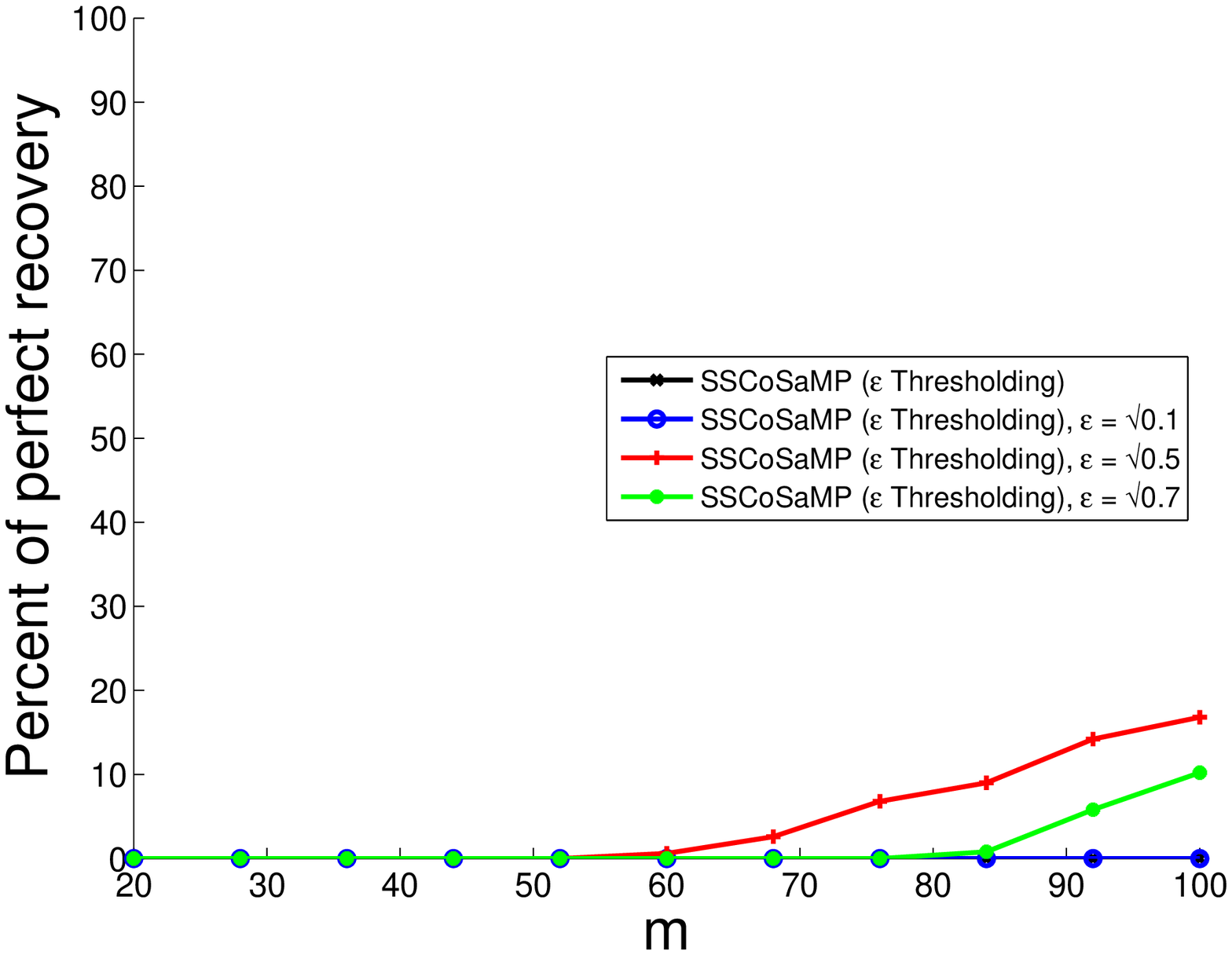}}
\end{minipage}
\caption{Recovery rate for SSCoSaMP ($\epsilon$-Thresholding)
with different values of $\epsilon$ for a random $m\times 1024$ Gaussian matrix $\matr{M}$
and a $4$ times overcomplete DFT matrix $\matr{D}$.  The signal is $8$-sparse and on the left the coefficients of the original signal are clustered whereas on the right they are separated.}
\label{fig:recovery rate2}
\end{figure}

\begin{figure}[htb]
\begin{minipage}[b]{.48\linewidth}
  \centering
  \centerline{\includegraphics[width=7.0cm]{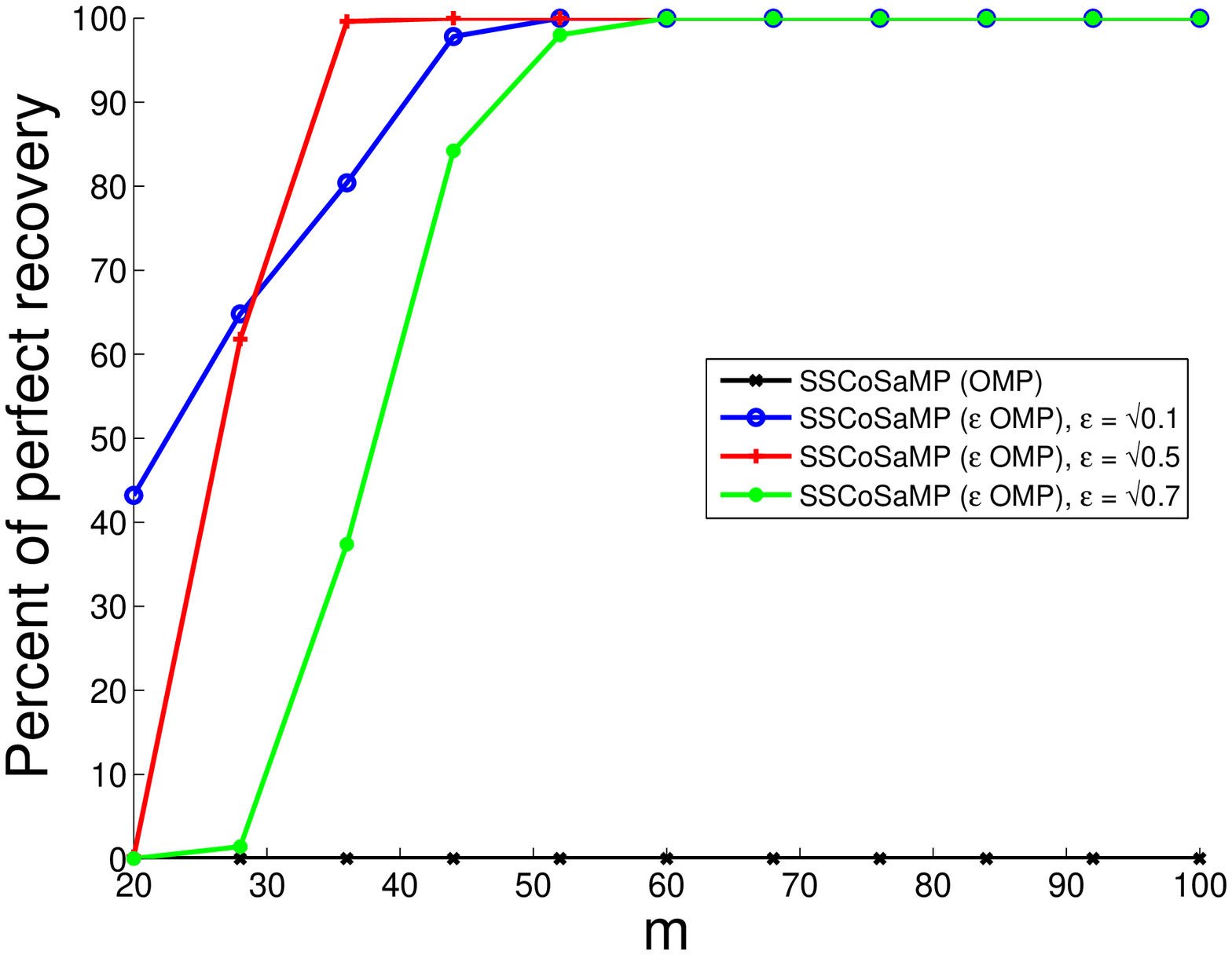}}
\end{minipage}
\hfill
\begin{minipage}[b]{.48\linewidth}
  \centering
  \centerline{\includegraphics[width=7.0cm]{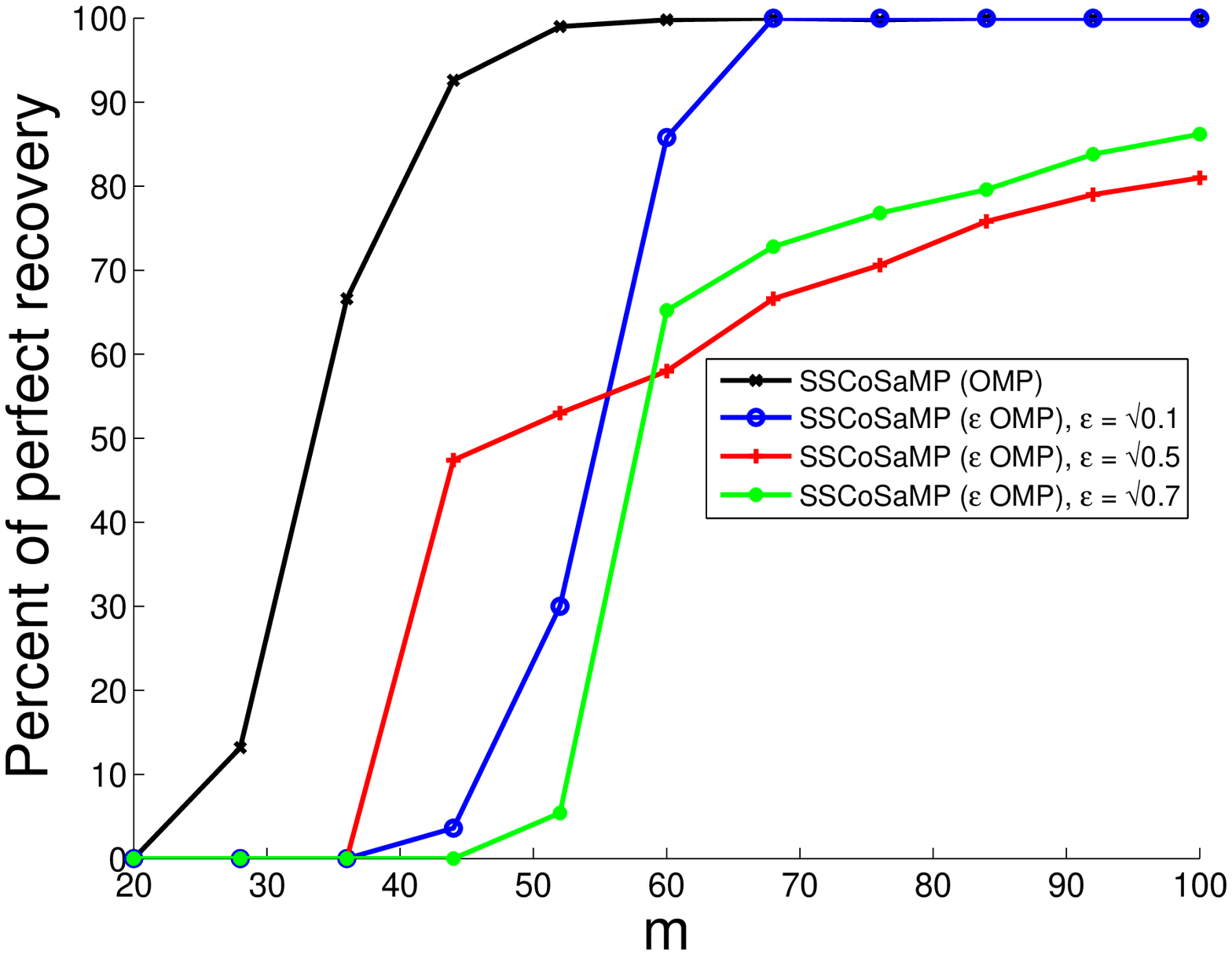}}
\end{minipage}
\caption{Recovery rate for SSCoSaMP ($\epsilon$-OMP)
with different values of $\epsilon$ for a random $m\times 1024$ Gaussian matrix $\matr{M}$
and a $4$ times overcomplete DFT matrix $\matr{D}$.  The signal is $8$-sparse and on the left the coefficients of the original signal are clustered whereas on the right they are separated.}
\label{fig:recovery rate3}
\end{figure}

The recovery results appear in Figures~\ref{fig:recovery rate}--\ref{fig:recovery rate3}.
 As seen from Figure~\ref{fig:recovery rate}, in the separated case SSCoSaMP-OMP works better for small values of $m$.  This is likely because it uses a smaller support set for which it is easier to satisfy the RIP condition.  As separated atoms are very uncorrelated it is likely that OMP will not be confused between them.
When the atoms are clustered, the high correlations take more effect and OMP is not able to recovery the right support because of the high coherence between close atoms in the cluster and around it. This is overcame by using $\epsilon$-OMP which uses larger support sets and thus resolves the confusion.  Note that even $\epsilon$-threshodling gets better recovery in this case, though it is a much simpler technique, and this shows that indeed the improvement is due to the extended support selection strategy.  As expected, using larger support estimates for the projection is more effective when the number of measurements $m$ is large.

We may say that the support extension step leads to a better recovery rate overall as it gets a good recovery on both the separated and clustered coefficient cases. In \cite{Davenport13Signal} it is shown that all the projection algorithms either perform well on the first case and very bad on the other or vice versa.
Using SSCoSaMP with $\epsilon$-OMP we have, at the cost of getting slightly inferior behavior in the separated case compared to SSCoSaMP with OMP, much improved behavior for the clustered case where the latter gets no recovery at all.

Figures~\ref{fig:recovery rate2}--\ref{fig:recovery rate3} demonstrate the sensitivity of the approximation algorithms to the choice of $\epsilon$ (note that $\epsilon=0$ reverts to the Thresholding/OMP algorithm).  While it is clear that $\epsilon$ cannot be too large (or far too many atoms will be included), the optimum choice of $\epsilon$ may not always be easy to identify since it depends on the dictionary $\matr{D}$.

\section{Discussion}
\label{sec:conc}

\subsection{Related Work}
Our work extends the work of Davenport, Needell, and Wakin~\cite{Davenport13Signal} who develop and analyze the Signal Space CoSaMP algorithm.  In that work, the $\matr{D}$-RIP is enforced, as well as access to projections which satisfy~\eqref{dnw_reqs}.  It is currently unknown whether there exist efficient projections which satisfy these requirements, even for well-behaved dictionaries like those that satisfy the RIP or have an incoherence property.  That being said, other results on signal space methods rely on such assumptions.  For example, a related work by Blumensath analyzes an algorithm which is a signal space extension of the Iterative Hard Thresholding (IHT) method~\cite{Blumensath11Sampling}.  The model in that work utilizes a union-of-subspaces model and also assumes the $\matr{D}$-RIP and projections with even stronger requirements than those in~\eqref{dnw_reqs}.

These types of projections also appear in \textit{model-based} compressive sensing, where such operators project onto a specified model set.  The model may describe structured sparsity patterns like tree-like or block sparsity, or may be a more general mode.  In this setting, signal recovery is performed by first reconstructing the coefficient vector, and then mapping to signal space.  When the dictionary $\matr{D}$ is an orthonormal basis, greedy methods have been adapted to structured sparsity models~\cite{baraniuk2010model}.  The assumptions, however, nearly require the product $\matr{AD}$ to satisfy the traditional RIP, and so extensions to non-orthonormal dictionaries serve to be difficult.  Although our work differs in its assumptions and domain model, model-based methods inspired the development of signal space CoSaMP~\cite{Davenport13Signal, davenport2012compressive}.

The importance of using and analyzing two separate projection schemes in sparse recovery is also discussed in an independent line of work by Hegde et.al.~\cite{hegdeapproximation}.  There, the authors call the two projections the ``head'' and ``tail'' projections, and analyze a variant of Iterative Hard Thresholding (IHT) for signal recovery under the Model-RIP, a generalization of the $\matr{D}$-RIP.  In fact, they show that without a projection satisfying essentially the second inequality of~\eqref{eq:C_optimal_proj}, conventional IHT will fail.

It would be also important to mention the relation of the $\epsilon$-OMP and $\epsilon$-thresholding algorithms (Algorithms~\ref{alg:OMP_eps} and \ref{alg:thresh_eps}) to the methods proposed in \cite{fannjiang2012CS, Duarte13Spectral}. The notion of excluding coherent atoms in the process of building the representation is used also within these works.
In particular, without the extension step, the $\epsilon$-OMP and $\epsilon$-threshodling techniques share a great similarity with the Band-Excluded OMP (BOMP) and Band-Excluded Matched Thresholding (BMT) methods in \cite{fannjiang2012CS} and the heuristic coherence-inhibiting sparse approximation strategy in \cite{Duarte13Spectral}. 
As we have seen in Section~\ref{sec:high_cor}, the use of the extension step deteriorates the performance in the case of separated coefficients as a larger support is processed and therefore the RIP conditions are harder to be satisfied. It is likely that using the techniques in \cite{fannjiang2012CS, Duarte13Spectral} would be better suited to deal with separated coefficient vectors. 

Finally, a related but significantly different vein of work also exists that studies signals from \textit{analysis} space rather than synthesis signal space.  Indeed, it was in this vein that the $\matr{D}$-RIP was first proposed and enforced for reconstruction~\cite{Candes11Compressed}.  In this setting, one requires that the analysis coefficients $\matr{D}^*\vect{x}$ are sparse or compressible, and reconstruction is performed in that domain.  Standard optimization based and greedy methods for compressed sensing have been extended and studied in this setting as well.  In particular, $\ell_1$-minimization~\cite{elad07Analysis, Candes11Compressed,Nam12Cosparse, Vaiter13Robust,Needell13Stable}, and greedy methods like CoSaMP and IHT have all been adapted to account for analysis (co)sparsity~\cite{Giryes13Greedy, Giryes12CoSaMP, Giryes11Iterative, Peleg12Performance}.

\subsection{Contributions and Directions}
In this work, we extend previous analysis of the Signal Space CoSaMP (SSCoSaMP) algorithm.  In signal space greedy algorithms for dictionaries which are not orthonormal, the computational bottleneck lies in the use of the approximate projections.  Here we extend the idea of a \textit{near-optimal projection}, and consider two possibly different near-optimal projections in the SSCoSaMP method.  Our new analysis enforces weaker assumptions on these projections, which hold when the dictionary $\matr{D}$ is incoherent or satisfies the RIP, unlike previous results whose assumptions do not hold in this setting.  Above, we discuss several important settings and describe algorithms that can be used for the approximate projections which satisfy our requirements for accurate signal recovery.  This includes even the case when the dictionary is highly coherent but each atom is only highly correlated with a small number of atoms, an important example in applications like super-resolution.

It remains an important and challenging open problem to develop approximate projection techniques which satisfy the assumptions of our main results even when the dictionary is highly coherent in an arbitrary fashion.  There are clearly limitations in this regard, as decoding from highly correlated atoms has fundamental theoretic boundaries.  It is unknown, however, how far these limits reach and for what applications accurate reconstruction is still possible.  An alternative of course is to develop greedy methods which do not require such projections, which we believe to be an equally challenging problem.


\section*{Acknowledgment} R. Giryes is grateful to the Azrieli Foundation for
the award of an Azrieli Fellowship.  D. Needell was partially supported by the
Simons Foundation Collaboration grant $\#274305$, the Alfred P. Sloan fellowship, 
and NSF Career grant $\#1348721$.  In addition, the authors thank the reviewers
of the manuscript for their suggestions which greatly improved the paper.

\appendix

\section{Proof of Lemma~\ref{lem:SSCoSaMP_xp_bound}}
\label{sec:SSCoSaMP_xp_bound_proof}

{\em Lemma~\ref{lem:SSCoSaMP_xp_bound}:}
If $\matr{M}$ has the $\matr{D}$-RIP with constants $\delta_{3\pscI k}, \delta_{(3\pscI +1)k}$, then
\begin{eqnarray*}
\norm{\vect{x}_p -\vect{x}}_2 &\le& \frac{1}{\sqrt{1-\delta_{(3\pscI +1)k}^2}}\norm{\Q_{\tilde{T}^t}(\vect{x}_p - \vect{x})}_2 +  \frac{\sqrt{1+\delta_{3\pscI k}}}{1-\delta_{(3\pscI +1)k}}\norm{\vect{e}}_2.
\end{eqnarray*}

{\em Proof:}
Since $\vect{x}_p \triangleq \matr{D}\alphabf_p$ is the minimizer of $\norm{\vect{y} - \vect{M}\tilde{\vect{x}}}_2$ with the constraints $\tilde{\vect{x}} = \matr{D}\tilde{\alphabf}$ and $\tilde{\alphabf}_{(\tilde{T}^t)^C} =0$, then
\begin{eqnarray}
\langle \matr{M} \vect{x}_p - \vect{y}, \matr{M}\vect{v} \rangle =0
\end{eqnarray}
for any vector $\vect{v} = \matr{D}\tilde\alphabf$ such that $\tilde\alphabf_{(\tilde{T}^t)^C} =0$.
Substituting $\vect{y} = \matr{M}\vect{x} + \vect{e}$ with simple arithmetics gives
\begin{eqnarray}
\label{eq:xp_x_property}
\langle \vect{x}_p - \vect{x}, \matr{M}^*\matr{M}\vect{v} \rangle = \langle \vect{e}, \matr{M}\vect{v} \rangle
\end{eqnarray}
where $\vect{v} = \matr{D}\tilde\alphabf$ and $\tilde\alphabf_{(\tilde{T}^t)^C} =0$.
Turning to look at $\norm{\P_{\tilde{T}^t}(\vect{x}_p - \vect{x})}_2^2$ and using \eqref{eq:xp_x_property} with $\vect{v} = \P_{\tilde{T}^t}(\vect{x}_p - \vect{x})$, we have
\begin{eqnarray}
\label{eq:Q_xp_x_norm}
&& \hspace{-0.3in} \norm{\P_{\tilde{T}^t}(\vect{x}_p - \vect{x})}_2^2 = \langle \vect{x}_p - \vect{x}, \P_{\tilde{T}^t}(\vect{x}_p - \vect{x}) \rangle \\
\nonumber && \hspace{-0.3in} = \langle \vect{x}_p - \vect{x}, (\matr{I}_d - \matr{M}^*\matr{M})\P_{\tilde{T}^t}(\vect{x}_p - \vect{x}) \rangle + \langle \vect{e}, \matr{M}\P_{\tilde{T}^t}(\vect{x}_p - \vect{x}) \rangle
\\ \nonumber && \hspace{-0.3in} \le \norm{ \vect{x}_p - \vect{x}}_2 \norm{\P_{\tilde{T}^t \cup T} (\matr{I}_d - \matr{M}^*\matr{M})\P_{\tilde{T}^t}}_2 \norm{\P_{\tilde{T}^t}(\vect{x}_p - \vect{x})}_2  \\ \nonumber && + \norm{ \vect{e}}_2\norm{ \matr{M}\P_{\tilde{T}^t}(\vect{x}_p - \vect{x})}_2
\\ \nonumber && \hspace{-0.3in} \le \delta_{(3\pscI +1)k}\norm{ \vect{x}_p - \vect{x}}_2 \norm{\P_{\tilde{T}^t}(\vect{x}_p - \vect{x})}_2  \\ \nonumber &&  + \norm{ \vect{e}}_2\sqrt{1+\delta_{3\pscI k}}\norm{\P_{\tilde{T}^t}(\vect{x}_p - \vect{x})}_2.
\end{eqnarray}
where the first inequality follows from the Cauchy-Schwartz inequality, the projection property that $\P_{\tilde{T}^t} = \P_{\tilde{T}^t}\P_{\tilde{T}^t}$ and the fact that $\vect{x}_p - \vect{x} = \P_{\tilde{T}^t \cup T}(\vect{x}_p - \vect{x})$.
The last inequality is due to the $\matr{D}$-RIP property, the fact that $|\tilde{T}^t| \le 3\pscI k$ and Corollary~\ref{cor:D_RIP_norm_diff}.
After simplification of \eqref{eq:Q_xp_x_norm} by $\norm{\P_{\tilde{T}^t}(\vect{x}_p - \vect{x})}_2$ we have
\begin{eqnarray}
\nonumber \norm{\P_{\tilde{T}^t}(\vect{x}_p - \vect{x})}_2 \le \delta_{(3\pscI +1)k}\norm{\vect{x}_p - \vect{x}}_2   + \sqrt{1+\delta_{3\pscI k}}\norm{\vect{e}}_2.
\end{eqnarray}
Utilizing the last inequality with the fact that $\norm{\vect{x}_p - \vect{x}}_2^2 = \norm{\Q_{\tilde{T}^t}(\vect{x}_p - \vect{x})}_2^2+ \norm{\P_{\tilde{T}^t}(\vect{x}_p - \vect{x})}_2^2$ gives
\begin{eqnarray}
&& \norm{\vect{x}_p - \vect{x}}_2^2  \le \norm{\Q_{\tilde{T}^t}(\vect{x}_p - \vect{x})}_2^2 + \left(\delta_{(3\pscI +1)k}\norm{\vect{x}_p - \vect{x}}_2   + \sqrt{1+\delta_{3\pscI k}}\norm{\vect{e}}_2 \right)^2.
\end{eqnarray}
The last equation is a second order polynomial of $\norm{\vect{x}_p - \vect{x}}_2$. Thus its larger root is an upper bound for it and this gives the inequality in \eqref{eq:SSCoSaMP_xp_bound}. For more details look at the derivation of (13) in \cite{foucart10Sparse}.
\hfill $\Box$ \bigskip

\section{Proof of Lemma~\ref{lem:SSCoSaMP_xt_bound1}}
\label{sec:SSCoSaMP_xt_bound1_proof}

{\em Lemma~\ref{lem:SSCoSaMP_xt_bound1}:}
Under the assumptions and notation of Theorem~\ref{thm:general_bound}, we have
\begin{eqnarray}
&& \hspace{-0.5in} \norm{\vect{x}^t -\vect{x}}_2 \le \rho_1\norm{\Q_{\tilde{T}^t}(\vect{x}_p - \vect{x})}_2+  \eta_1\norm{\vect{e}}_2
\end{eqnarray}

{\em Proof:}

We start with the following observation
\begin{eqnarray}
\label{eq:xt_x_diff_norm}
&& \hspace{-0.3in} \norm{\vect{x} -\x^t}_2 = \norm{\vect{x} - \vect{x}_p +  \vect{x}_p-\x^t}_2
 \le \norm{\vect{x}- \vect{x}_p}_2 + \norm{\x^t -  \vect{x}_p}_2,
\end{eqnarray}
where the last step is due to the triangle inequality. 
Using \eqref{eq:C_optimal_ineq} with the fact that $\x^t = \P_{\CFI_{B,\pscI k}( \vect{x}_p)} \vect{x}_p$ we have
\begin{eqnarray}
\label{eq:xt_x_diff_norm_coeff1}
\norm{\x^t -  \vect{x}_p}_2^2 \le \CI_k\norm{\x - \vect{x}_p}_2^2.
\end{eqnarray}
Plugging \eqref{eq:xt_x_diff_norm_coeff1} in \eqref{eq:xt_x_diff_norm} leads to 
\begin{eqnarray}
\label{eq:xt_x_diff_norm_all}
\norm{\x -\x^t}_2  &\le & (1 + \sqrt{C_k})\norm{\vect{x} - \vect{x}_p}_2 \\ \nonumber &\le & \frac{1 + \sqrt{C_k}}{\sqrt{1-\delta_{(3\psc +1)k}^2}}\norm{\Q_{\tilde{T}^t}( \vect{x}_p - \vect{x})}_2 
+  \frac{\left(1 + \sqrt{C_k}\right)\sqrt{1+\delta_{3\pscI k}}}{1-\delta_{(3\psc +1)k}}\norm{\P_{T_{\vect{e}}}\matr{M}^*\vect{e}}_2,
\end{eqnarray}
where for the last inequality we use Lemma~\ref{lem:SSCoSaMP_xp_bound}.

\hfill $\Box$ \bigskip

\section{Proof of Lemma~\ref{lem:SSCoSaMP_Pxp_bound}}
\label{sec:SSCoSaMP_Pxp_bound_proof}

{\em Lemma~\ref{lem:SSCoSaMP_Pxp_bound}:}
Under the assumptions and notation of Theorem~\ref{thm:general_bound}, we have
\begin{eqnarray}
&& \hspace{-0.3in}\norm{\Q_{\tilde{T}^t}(\vect{x}_p - \vect{x})}_2 \le 
 \eta_2\norm{\vect{e}}_2
 + \rho_2\norm{\vect{x} - \vect{x}^{t-1}}_2.
\end{eqnarray}

{\em Proof:}
Looking at the step of finding new support elements one can observe that $\P_{T_{\Delta}}$ is a near optimal projection operator for $\matr{M}^*\vect{y}^{t-1}_r = \matr{M}^*(\vect{y} - \matr{M}\vect{x}^{t-1})$.
Noticing that $T_{\Delta} \subseteq \tilde{T}^t$ and then using \eqref{eq:C_optimal_ineq_up} with $\P_{T^{t-1} \cup T}$ gives
\begin{eqnarray}
\label{eq:SSCoSaMP_PMy_Mx_ineq}
&& \hspace{-0.3in} \norm{\P_{\tilde{T}^t}\matr{M}^*(\vect{y} - \matr{M}\vect{x}^{t-1})}_2^2 \\ \nonumber && ~~~~~~~~~~~~~~~~~~~~~~ \ge \norm{\P_{{T}_{\Delta}}\matr{M}^*(\vect{y} - \matr{M}\vect{x}^{t-1})}_2^2
\\ \nonumber && ~~~~~~~~~~~~~~~~~~~~~~\ge \CII_{2k}\norm{\P_{T^{t-1} \cup T}\matr{M}^*(\vect{y} - \matr{M}\vect{x}^{t-1})}_2^2.
\end{eqnarray}

We start by bounding the lhs of \eqref{eq:SSCoSaMP_PMy_Mx_ineq} from above. Using Proposition~\ref{prop:norm2_ineq}
with $\gamma_1>0$ and $\alpha>0$ we have
\begin{eqnarray}
\label{eq:SSCoSaMP_PMy_Mx_ineq_lhs}
&& \hspace{-0.3in} \norm{\P_{\tilde{T}^t}\matr{M}^*(\vect{y} - \matr{M}\vect{x}^{t-1})}_2^2 \le
(1+\frac{1}{\gamma_1})\norm{\P_{\tilde{T}^t}\matr{M}^*\vect{e}}_2^2
\\ \nonumber && + (1+\gamma_1)\norm{\P_{\tilde{T}^t}\matr{M}^*\matr{M}(\vect{x} - \vect{x}^{t-1})}_2^2
\\ \nonumber && \hspace{-0.3in} \le \frac{1+\gamma_1}{\gamma_1}\norm{\P_{\tilde{T}^t}\matr{M}^*\vect{e}}_2^2
 + (1+\alpha)(1+\gamma_1)\norm{\P_{\tilde{T}^t}(\vect{x} - \vect{x}^{t-1})}_2^2  \\ \nonumber &&
+(1+\frac{1}{\alpha})(1+\gamma_1)\norm{\P_{\tilde{T}^t}(\matr{I}_d-\matr{M}^*\matr{M})(\vect{x} - \vect{x}^{t-1})}_2^2
\\ \nonumber && \hspace{-0.3in} \le  \frac{(1+\gamma_1)(1+\delta_{3\pscI k})}{\gamma_1}\norm{\vect{e}}_2^2
 \\ \nonumber && - (1+\alpha)(1+\gamma_1)\norm{\Q_{\tilde{T}^t}(\vect{x} - \vect{x}^{t-1})}_2^2
 \\ \nonumber &&  +\left(1+\alpha +\delta_{(3\pscI+1)k} + \frac{\delta_{(3\pscI +1)k}}{\alpha}\right)(1+\gamma_1)\norm{\vect{x} - \vect{x}^{t-1}}_2^2,
\end{eqnarray}
where the last inequality is due to Corollary~\ref{cor:MP_RIP_norm} and \eqref{eq:D_RIP_norm_diff}.

We continue with bounding the rhs of \eqref{eq:SSCoSaMP_PMy_Mx_ineq} from below.
For the first element we use Proposition~\ref{prop:norm2_ineq} with constants $\gamma_2>0$ and $\beta >0$, and \eqref{eq:D_RIP_norm_diff} to achieve
\begin{eqnarray}
\label{eq:SSCoSaMP_PMy_Mx_ineq_rhs1}
&& \hspace{-0.3in} \norm{\P_{{T^{t-1} \cup T}}\matr{M}^*(\vect{y} - \matr{M}\vect{x}^{t-1})}_2^2
\\ \nonumber && \hspace{-0.3in} \ge \frac{1}{1+\gamma_2}\norm{\P_{T^{t-1} \cup T}\matr{M}^*\matr{M}(\vect{x} - \vect{x}^t)}_2^2
-\frac{1}{\gamma_2}\norm{\P_{T^{t-1} \cup T}\matr{M}^*\vect{e}}_2^2
\\ \nonumber && \hspace{-0.3in} \ge \frac{1}{1+\beta}\frac{1}{1+\gamma_2}\norm{\vect{x} - \vect{x}^{t-1}}_2^2
-\frac{1}{\gamma_2}\norm{\P_{T^{t-1} \cup T}\matr{M}^*\vect{e}}_2^2
\\ \nonumber && -\frac{1}{\beta}\frac{1}{1+\gamma_2}\norm{\P_{T^{t-1} \cup T}(\matr{M}^*\matr{M} - \matr{I}_d)(\vect{x} - \vect{x}^{t-1})}_2^2
\\ \nonumber && \hspace{-0.3in} \ge (\frac{1}{1+\beta}-\frac{\delta_{(\pscI+1)k}}{\beta})\frac{1}{1+\gamma_2}\norm{\vect{x} - \vect{x}^{t-1}}_2^2
-\frac{(1+\delta_{(\pscI+1)k})}{\gamma_2}\norm{\vect{e}}_2^2.
\end{eqnarray}

By combining \eqref{eq:SSCoSaMP_PMy_Mx_ineq_lhs} and \eqref{eq:SSCoSaMP_PMy_Mx_ineq_rhs1} with
\eqref{eq:SSCoSaMP_PMy_Mx_ineq} we have
\begin{eqnarray}
&& (1+\alpha)(1+\gamma_1)\norm{\Q_{\tilde{T}^t}(\vect{x} - \vect{x}^{t-1})}_2^2  \\ \nonumber &&
 \le  \frac{(1+\gamma_1)(1+\delta_{3\pscI k})}{\gamma_1}\norm{\vect{e}}_2^2
 +\CII_{2k}\frac{(1+\delta_{(1+\pscI)k})}{\gamma_2}\norm{\vect{e}}_2^2
 \\ \nonumber &&
 + \left(1+\alpha +\delta_{(3\pscI+1)k} + \frac{\delta_{(3\pscI+1)k}}{\alpha}\right)(1+\gamma_1)\norm{\vect{x} - \vect{x}^{t-1}}_2^2
\\ \nonumber &&
- \CII_{2k}(\frac{1}{1+\beta}-\frac{\delta_{(1+\pscI)k}}{\beta})\frac{1}{1+\gamma_2}\norm{\vect{x} - \vect{x}^{t-1}}_2^2.
\end{eqnarray}
Division of both sides by $(1+\alpha)(1+\gamma_1)$ yields
\begin{eqnarray}
&& \hspace{-0.3in}\norm{\Q_{\tilde{T}^t}(\vect{x} - \vect{x}^{t-1})}_2^2 \le \\ \nonumber &&
 \hspace{-0.2in}   \bigg(\frac{1+\delta_{3\pscI k}}{\gamma_1(1+\alpha)}
 +\frac{(1+\delta_{(\pscI+1)k})\CII_{2k}}{\gamma_2(1+\alpha)(1+\gamma_1)}
\bigg)\norm{\vect{e}}_2^2
 \\ \nonumber && \hspace{-0.3in}
 + \bigg( 1 +\frac{\delta_{(3\pscI+1)k}}{\alpha}
 \\ \nonumber && \hspace{-0.2in} - \frac{\CII_{2k}}{(1+\alpha)(1+\gamma_1)(1+\gamma_2)}(\frac{1}{1+\beta}-\frac{\delta_{(\pscI+1)k}}{\beta}) \bigg)\norm{\vect{x} - \vect{x}^{t-1}}_2^2.
\end{eqnarray}
Substituting $\beta = \frac{\sqrt{\delta_{(\pscI+1)k}}}{1- \sqrt{\delta_{(\pscI+1)k}}}$ gives
\begin{eqnarray}
&& \hspace{-0.3in}\norm{\Q_{\tilde{T}^t}(\vect{x} - \vect{x}^{t-1})}_2^2 \le \\ \nonumber &&
 \hspace{-0.2in}   \bigg(\frac{1+\delta_{3\pscI k}}{\gamma_1(1+\alpha)}
 +\frac{(1+\delta_{(\pscI+1)k})\CII_{2k}}{\gamma_2(1+\alpha)(1+\gamma_1)}\bigg)\norm{\vect{e}}_2^2
 \\ \nonumber && \hspace{-0.3in}
 + \bigg(1 +\frac{\delta_{(3\pscI+1)k}}{\alpha}
 \\ \nonumber && \hspace{-0.2in}- \frac{\CII_{2k}}{(1+\alpha)(1+\gamma_1)(1+\gamma_2)}\left(1-\sqrt{\delta_{(\pscI+1)k}}\right)^2\bigg)\norm{\vect{x} - \vect{x}^{t-1}}_2^2,
\end{eqnarray}
Using $  {\alpha = \frac{\sqrt{\delta_{(3\pscI+1)k}}}{\sqrt{\frac{\CII_{2k}}{(1+\gamma_1)(1+\gamma_2)}}\left(1-\sqrt{\delta_{(\pscI+1)k}}\right)- \sqrt{\delta_{(3\pscI+1)k}}}}$ yields
\begin{eqnarray}
&& \hspace{-0.6in}\norm{\Q_{\tilde{\Lambda}^t}(\vect{x} - \vect{x}^{t-1})}_2^2 \le \\ \nonumber &&
 \hspace{-0.6in}   \bigg(\frac{1+\delta_{3\pscI k}}{\gamma_1(1+\alpha)}
 +\frac{(1+\delta_{(\pscI+1)k})\CII_{2k}}{\gamma_2(1+\alpha)(1+\gamma_1)}\bigg)\norm{\vect{e}}_2^2
 \\ \nonumber && \hspace{-0.6in}
 + \bigg(-\bigg( \sqrt{\delta_{(3\pscI+1)k}}-\sqrt{\frac{\CII_{2k}}{(1+\gamma_1)(1+\gamma_2)}}\left(1-\sqrt{\delta_{(\pscI+1)k}}\right)
 \bigg)^2
 \\ \nonumber && \hspace{-0.6in} +1\bigg)\norm{\vect{x} - \vect{x}^{t-1}}_2^2,
\end{eqnarray}
The values of $\gamma_1, \gamma_2$ give a tradeoff between the convergence rate and the size of the noise coefficient.
For smaller values we get better convergence rate but higher amplification of the noise.
We make no optimization on them and choose them to be $\gamma_1 = \gamma_2 = \gamma$ where $\gamma$ is an arbitrary number greater than $0$.
Thus we have
\begin{eqnarray}
&& \hspace{-0.6in}\norm{\Q_{\tilde{T}^t}(\vect{x} - \vect{x}^{t-1})}_2^2 \le \\ \nonumber &&
 \hspace{-0.6in}   \bigg(\frac{1+\delta_{3\pscI k}}{\gamma(1+\alpha)}
 +\frac{(1+\delta_{(\pscI+1)k})\CII_{2k}}{\gamma(1+\alpha)(1+\gamma)}\bigg)\norm{\vect{e}}_2^2
 \\ \nonumber && \hspace{-0.6in}
 + \bigg(-\bigg( \sqrt{\delta_{(3\pscI+1)k}}-\frac{\sqrt{\CII_{2k}}}{1+\gamma}\left(1-\sqrt{\delta_{(\pscI+1)k}}\right)
 \bigg)^2
 \\ \nonumber && \hspace{-0.6in} +1\bigg)\norm{\vect{x} - \vect{x}^{t-1}}_2^2,
\end{eqnarray}
Using the triangle inequality and the fact that $\Q_{\tilde{T}^t}\vect{x}_p = \Q_{\tilde{T}^t}\vect{x}^{t-1} =0$ gives the desired result.

\hfill $\Box$ \bigskip

\bibliographystyle{elsarticle-num}
\bibliography{SSCoSaMP_abb,IEEEabrv}

\end{document}